\documentclass[pdflatex]{svproc}
\usepackage{lmodern}

\usepackage[bookmarks,bookmarksopen=true,linkcolor=cyan,urlcolor=blue,citecolor=magenta]{hyperref}
\usepackage{fullpage}
\usepackage{subcaption}
\usepackage{mathtools,etoolbox,faktor}

\usepackage{booktabs,tabularx,float}
\usepackage{longtable}
\usepackage{todonotes}

\usepackage{xcolor}
\usepackage{graphicx,xparse,dirtytalk,soul}
\usepackage{comment}
\usepackage{xifthen}
\usepackage{ulem}
\usepackage{listings}

\usepackage{tikz,pgf}
\usetikzlibrary{arrows,patterns,positioning,shapes.multipart, fit,backgrounds,calc,decorations.pathreplacing, graphs, graphs.standard, shapes,decorations,angles}
\pgfdeclaredecoration{sl}{initial}{
  \state{initial}[width=\pgfdecoratedpathlength-1sp]{
     \pgfmoveto{\pgfpointorigin}
  }
  \state{final}{
     \pgflineto{\pgfpointorigin}
    }
}

\usepackage{thmtools}
\usepackage[english, ruled, onelanguage, linesnumbered, vlined]{algorithm2e}

\usepackage{import}
\usepackage[usestackEOL]{stackengine}
\usepackage{enumitem}

\usepackage{dsfont}
\usepackage{framed}
\usepackage[bb=boondox]{mathalfa}
\usepackage{fontawesome5}

\usepackage{tcolorbox}

\usepackage{rotating}

\SetKwComment{com}{$\triangleright$\ }{}
\usepackage{centernot}

\usepackage{pdfpages}
\usepackage{csquotes}

\usepackage{pgf,pgffor}

\usepackage{chngcntr}
\counterwithin{table}{section}

\usepackage[numbers]{natbib} 

\spnewtheorem{claimnum}{Claim}{\bfseries}{\itshape}

\usepackage[nameinlink,noabbrev, capitalize]{cleveref}
\AtBeginEnvironment{appendices}{\crefalias{section}{appendix}}
\crefname{subsection}{Subsection}{Subsections}
\crefname{subsubsection}{Subsubsection}{Subsubsections}
\crefname{proposition}{Proposition}{Propositions}
\crefname{lemma}{Lemma}{Lemmas}
\crefname{corollary}{Corollary}{Corollaries}
\crefname{claimnum}{Claim}{Claims}
\makeatletter
\newcommand{\labelx}[1]{
    \relax
    \ifmmode
        \label{#1} 
    \else 
        \ifnum\pdfstrcmp{\@currenvir}{document}=0
            \label{#1}
        \else
            \label[\@currenvir]{#1}
        \fi
    \fi
}
\makeatother

\newcommand{\bb}[1]{%
  \renewcommand*{\do}[1]{
    \expandafter\newcommand\csname ##1\endcsname{\ensuremath{\mathbb{##1}}}%
  }%
  \docsvlist{#1}
}
\newcommand{\opname}[1]{%
  \renewcommand*{\do}[1]{
    \expandafter\newcommand\csname ##1\endcsname{\ensuremath{\operatornamewithlimits{##1}}}%
  }%
  \docsvlist{#1}
}

\bb{R,N,Z,Q,C, e, p}
\opname{argmin, argmax, arglatest, latest}

\newcommand{\expc}{\textsc{exptime-complete}}
\newcommand{\uv}[2]{#1#2}
\newcommand{\duv}[2]{(#1, #2)}

\newcommand{\dg}{\mathcal{D}}
\newcommand{\tg}{\mathcal{G}}

\newcommand{\Sdo}{G_{\scriptscriptstyle \mathrm{min}}}
\newcommand{\Sup}{G_{\scriptscriptstyle \mathrm{max}}}

\newcommand{\domnum}{\gamma}

\newcommand{\pet}{P_{e}} 

\newcommand{\rad}{r}

\newcommand{\npg}[3]{\Gamma_{#1}\left(#2, #3\right)}

\newcommand{\ctgmax}[1]{c_{\text{\tiny{\faClock[regular]}}}(#1)}

\newcommand{\peri}[1]{\left(#1\right)^{*}}

\newcommand{\conf}{C}

\newcommand{\neig}[3]{N_{#1}\left[#2,#3\right]}

\newcommand{\cs}{\sigma_c}

\newcommand{\copnum}[1][]{%
  \ifthenelse{\isempty{#1}}%
  {cop number}%
  {cop numbers}%
}
\newcommand{\Copnum}[1][]{%
  \ifthenelse{\isempty{#1}}%
  {Cop number}%
  {Cop numbers}%
}
\newcommand{\pnei}[1][]{%
    \ifthenelse{\isempty{#1}}%
    {$p$-outneighbour}%
    {$p$-outneighbours}%
}

\newcommand{\tw}{\mathrm{tw}}
\newcommand{\tb}{\mathcal{B}}

\DeclarePairedDelimiterX{\set}[1]{\{}{\}}{\setargs{#1}}
\NewDocumentCommand{\setargs}{>{\SplitArgument{1}{;}}m}
{\setargsaux#1}
\NewDocumentCommand{\setargsaux}{mm}
{\IfNoValueTF{#2}{#1} {#1\,\delimsize|\,\mathopen{}#2}}

\newrobustcmd*{\myline}[1]{\tikz{\filldraw[draw=#1, line width = 3pt] (0,0) -- (0.5cm, 0);}}
\newrobustcmd*{\myfline}[1]{\tikz{\filldraw[draw=#1, line width = 0.2mm] (0,0) -- (0.5cm, 0);}}
\newrobustcmd*{\mycircle}[1]{\tikz{
  \node[draw=darkred, double, circle] (0,0) {};
  }
}
\newrobustcmd*{\mydot}[1]{\tikz{
  \node[fill=#1, circle] (0,0) {};
  }
}

\DeclarePairedDelimiter{\abs}{\lvert}{\rvert}
\makeatletter
\let\oldabs\abs
\def\abs{\@ifstar{\oldabs}{\oldabs*}}

\DeclarePairedDelimiter{\fl}{\lfloor}{\rfloor}
\makeatletter
\let\oldfl\fl
\def\fl{\@ifstar{\oldfl}{\oldfl*}}

\DeclarePairedDelimiter{\cl}{\lceil}{\rceil}
\makeatletter
\let\oldcl\cl
\def\cl{\@ifstar{\oldcl}{\oldcl*}}

\patchcmd{\abstract}{\titlepage}{\clearpage}{}{}

\makeatletter
\pgfmathdeclarefunction{alpha}{1}{%
  \pgfmathint@{#1}%
  \edef\pgfmathresult{\pgffor@alpha{\pgfmathresult}}%
}

\begin{document}

\title{Cop numbers of periodic graphs}
\titlerunning{Cop numbers of periodic graphs}

\author{
	Jean-Lou De Carufel\inst{1} \and
	Paola Flocchini\inst{2} \and
	Nicola Santoro\inst{3} \and
	Fr\'ed\'eric Simard\inst{4}
}
\authorrunning{Jean-Lou De Carufel et al.}
\institute{
	University of Ottawa, Ottawa, Ontario, Canada,\\ \email{jdecaruf@uottawa.ca} \and
	University of Ottawa, Ottawa, Ontario, Canada,\\ \email{paola.flocchini@uottawa.ca} \and
	Carleton University, Ottawa, Ontario, Canada,\\ \email{santoro@scs.carleton.ca} \and
	University of Ottawa, Ottawa, Ontario, Canada,\\ \email{fsima063@uottawa.ca}
}


\maketitle

\begin{abstract}
A \emph{periodic graph} ${\cal G}=(G_0, G_1, G_2, \dots)$ with period $p$ is an infinite periodic sequence of graphs $G_i = G_{i + p} = (V,E_i)$, where $i \geq 0$. The graph $G=(V,\cup_i E_i)$ is called the footprint of $\tg$. Recently, the arena where the  Cops and Robber game is  played has been extended  from a graph to a periodic graph; in this case, the \emph{\copnum{}} is also the minimum number of cops sufficient for capturing the robber. We study the connections and distinctions between the  \copnum{}  $c({\cal G})$ of a periodic graph ${\cal G}$ and the \copnum{} $c(G)$ of its footprint $G$ and establish several facts. For instance, we show that the smallest periodic graph with $c({\cal G}) = 3$  has at most $8$ nodes; in contrast, the smallest graph $G$ with $c(G) = 3$  has $10$ nodes. 

We push this investigation by generating multiple examples showing how the \copnum{s} of a periodic graph $\tg$, the subgraphs $G_i$ and its footprint $G$ can be loosely tied.

Based on these results, we derive upper bounds on the \copnum{} of a periodic graph from properties of its footprint such as its treewidth. 

\end{abstract}

\section{Introduction}

The game of Cops and Robber is a pursuit-evasion game played in turns, originally on a finite undirected graph $G$, between $k\geq 1$ cops and a single robber. There is perfect information in the sense that all players know $G$ and the locations of the other players at every turn. Initially, first the cops, then the robber, choose their positions on $G$. Then, in every turn each cop first moves to a neighbouring vertex or stays still, then the robber moves to a neighbouring vertex or stays still. The game ends and the $k$ cops win if they ever step on the node occupied by the robber. The robber wins by forever evading capture. 

This game was first described by Nowakowski and Winkler \cite{Nowakowski1983} and independently by Quilliot \cite{Qu78} for $k = 1$. Later, Aigner and Fromme \cite{Aigner1984} extended the game to general values of $k$. The smallest integer $k\geq 1$ for which $k$ cops can always capture the robber on $G$ is called the \copnum{} of $G$, denoted as $c(G)$. Determining whether $c(G)\leq k$ for an input of $(G,k)$ is \expc{} in general \cite{Kinnersley2015}. Nevertheless, computing $c(G)$, finding good upper and lower bounds on this number and finding classes of graphs for which $c(G)$ is bounded (by a constant) are some of the main research objectives in this field. 

Bonato and Nowakowski \cite{Bonato2011f} summarize many results on the Cops and Robber game in their seminal book. 

\subsection{Temporal Graphs}

A common feature of such disparate networks as wireless ad-hoc networks and social networks is that they are continuously subject to change. Such networks have spurred recent research into temporal, or time-varying, graphs and multiple models have been suggested to describe them (see for example \cite{Casteigts2013,Holme2012,WeZF15}).

In order to simplify the investigations on time-varying graphs, one generally imposes assumptions on time. One common assumption is to consider time as \emph{discrete}, like in the Cops and Robber games where the players move in turns. This often comes with the assumption that the number of nodes does not grow infinitely. In this case, a temporal graph is often described as an infinite sequence $\tg = (G_0,G_1,\dots)$ of subgraphs $G_i = (V,E_i)$ of a common graph $G = (V,\cup_i E_i)$, a model that was originally described in \cite{Fe04,HarG97}. Here, $G_t$ is called the \emph{snapshot} of $G$ at time $t$ while $G$ is the footprint.

It is also customary to apply assumptions on the variability of the changes. Some assumptions relate to the \say{connectivity} of $G$ and each $G_t$. The strongest condition, \emph{$1$-interval connectivity} (e.g., \cite{IlKW14,KuhLyO10,OdW05}), requires that each subgraph is connected. On the other hand, many weak conditions exist (see, e.g., \cite{Casteigts2013}), for example requiring only that the sequence be \emph{connected over time} (\cite{CaFMS14,GoSOKM21}). One can also wish to control the frequency of appearance of the edges in the sequence. The most relevant such assumption is \emph{periodicity}: there exists a positive integer $p$ such that $G_i = G_{i+p}$ for all $i\in \Z^+$ (e.g., \cite{flocchini2013exploration,IlW11,JatYG14}), so that $\tg$ can be written as a finite sequence $\tg = \peri{G_0,\dots,G_{p-1}}$. We refer to such temporal graphs as \emph{periodic graphs}.

There are multiple problems about \emph{mobile entities} in which the agents operate on temporal graphs under different conditions, such as \emph{graph exploration, dispersion} and \emph{gathering} (e.g., \cite{AgMMS18,BouDP17,DiLDFS20,DiLFPPSV17,erlebach2021temporal,ErS18,GoSOKM21,GoSOKM18} and \cite{DiL19} for a recent survey). Only recently was the Cops and Robber game, which can be considered as a really specific model of mobile agents, studied on (periodic) temporal graphs.

Extending the game of Cops and Robber to periodic graphs is straightforward because it is generally played on a finite structure and a periodic graph $\tg = \peri{G_0,\dots,G_{p-1}}$ can be thought of as a finite sequence of finite graphs. The rules are easily extended. All players know $\tg$ as well as the positions of the other players. Initially, the cops first choose a set of nodes in $V$ to occupy, then the robber chooses a vertex of $V$. Then, starting in $G_0$, first the cops, then the robber move to a node that is adjacent to their current positions. After both players (cops and robber) have moved in $G_t$, they start their next turn occupying their nodes in $G_{t+1}$ where they will play their next move. The game ends and the cops win if and only if one cop can move on the node occupied by the robber in some snapshot. The robber wins otherwise.

\subsection{Existing Results On Cops and Robber Games}

The literature on Cops and Robber games is vast, so we only highlight some interesting contributions that we use here. 
First, Clarke \cite{clarke2002constrained} showed that outerplanar graphs have \copnum{} at most $2$. This result came after Aigner and Fromme \cite{Aigner1984} had showed that planar graphs have \copnum{} at most $3$. 
Maamoun and Meyniel \cite{maamoun1987game} showed that $c(Q_k) = \cl{\frac{k + 1}{2}}$, where $Q_k$ is the hypercube of dimension $k$. One important concept is that of the \emph{retract}. A \emph{retraction} $h : G\to H$ is a homomorphism from $G$ to one of its subgraphs $H$ such that $h(H) = H$. The graph $H$ itself is called a retract of $G$. Retracts are used for example to show that graphs with \copnum{} one (also called \emph{copwin}) are exactly the \emph{dismantlable} graphs \cite{Nowakowski1983}. In particular, Berarducci and Intrigila \cite{berarducci1993cop} showed that whenever $H$ is a retract of $G$, $c(H)\leq c(G)$ so the \copnum{} is upper bounded under taking retracts (this is not the case for general subgraphs). This last result was used by Baird et al. \cite{baird2013minimum} and Turcotte and Yvon \cite{turcotte20214} to prove the size of the smallest graphs with \copnum{} $3$ and $4$ (respectively). Finally, Joret et al. \cite{joret2008cops} noticed an interesting connection between tree decompositions and Cops and Robber game and showed that $c(G)\leq \tw(G)/2 + 1$.

It is common to study the game of Cops and Robber either in different settings or with different rules. For example, one can let the robber move faster (see for example \cite{mehrabian2011cops}) or play on a directed graph. Simply increasing the speed of the robber leads to fascinating results. One can think of Fomin et al.'s \cite{fomin2010pursuing} result that the $n\times n$ grid has \copnum{} $\Omega(\sqrt{\log{n}})$ when the robber has speed two. 

The game on a directed graph has been a subject of much research because of how difficult it is to analyze. Khatri et al. \cite{Khatri2018} present many interesting results. Loh and Oh \cite{loh2017cops} exhibited a planar directed graph with \copnum{} $4$, which shows that Aigner and Fromme's upper bound on planar graphs does not hold in the directed case. 

In order to help analyze all variants of Cops and Robber games with a single framework, Bonato and MacGillivray \cite{BoM17} and, separately, Simard et al. \cite{simard2021general} described models and algorithms to compute the \copnum{} of any Cops and Robber games, under some conditions. While Bonato and MacGillivray focus on deterministic games, Simard et al.'s model can handle probabilistic outcomes.

\subsection{Existing Results On Cops and Robbers Games on Temporal Graphs}

The game of Cops and Robber was first brought to the context of periodic graphs by \citet{ErS20}. They studied this game with a single cop and showed an algorithm to determine if a periodic graph has \copnum{} one or more. They also have some results about periodic graphs where the footprint is a cycle. Their algorithm works by transforming the game into a reachability game \cite{Berwanger2013GraphGW}. Bavel et al. \cite{Balev2020} similarly conceived algorithms to determine whether a periodic graph has \copnum{} one or not. Interestingly, they also studied an \say{online} version of the game in which the sequence of graphs is not known a priori. This game has imperfect information, which makes the problem of computing the \copnum{} much harder. Other authors had previously studied restrictions on information available to the cops in the game played on undirected graphs, such as Clarke \cite{Clarke2001}. Nevertheless, it is much more common to assume perfect information.
Morawietz et al. \cite{morawietz2020,morawietz2021timecop} studied different parameterized complexity problems related to computing the \copnum{} of a periodic graph. They showed that if the periodic graph is not given explicitly as a sequence of graphs, but only as a boolean \say{edge-presence function} (see \cite{Casteigts2013}), then determining if a periodic graph has \copnum{} one or not is \textsc{np-hard}. This question was raised in \cite{ErS20}. The results and algorithmic questions presented in \cite{ErS20,morawietz2020,morawietz2021timecop} have been merged and discussed  in a more complete form in  \citet{erlebach2024}.

\citet{CaFlSaSi23} gave a characterization of periodic graphs with \copnum{} one as well as an algorithm to determine if a periodic graph has \copnum{} one that is more efficient than the one of \citet{ErS20}. A thorough review of the literature on Cops and Robber games on static and temporal graphs as well as the  results of 
\cite{CaFlSaSi23} also appear in \cite{simard2024temporal}, along with other novel results.

\subsection{Contributions}

Notably absent from the state of current research on the Cops and Robber game on periodic graphs are deep analytical results about periodic graphs with varying \copnum{s}. With the aim of understanding the relationship between the static and temporal settings in terms of this classical parameter, we start by focusing on the differences between the \copnum{} of a periodic graph and the \copnum{s} of its constituent static graphs. Our results show that the temporal dimension introduces huge differences from the static setting, and we discover some properties of graphs that help control those variations. 

We investigate the maximum \copnum{} of any periodic graph with footprint $G$, written as $\ctgmax{G}$. From this quantity we seek to understand how the footprint constrains the \copnum{} of a periodic graph. One takeaway from our investigation of $\ctgmax{G}$ is that copwin strategies on periodic graphs, when only $G$ is known, need to be resilient to change. This setting is akin to planning under uncertainty, when failures can occur in the graph $G$, and this is easier when $G$ has good separation properties. 

We contrast Baird et al's result \cite{baird2013minimum} on the minimum order of a $3$-copwin graph by exhibiting a smaller (in the order of its footprint) periodic graph with the same \copnum{} in \cref{prop:hypercubeQ3_3copwin}. We show that no value of $c(G), c(G_0),\dots,c(G_{p-1})$ can, in general, be used as either lower or upper bound on $c(\tg)$ by presenting two counterexamples (\cref{conj:112_periodicgraph} and \cref{thm:up_low_bounds_impos}). We complete this presentation by filling \cref{tab:summ_results} that summarizes many examples we present and their different \copnum{s}. This table serves to highlight the counterintuitive nature of periodic graphs and the difficulty of deriving $c(\tg)$ from $c(G), c(G_0),\dots,c(G_{p-1})$. Those examples are presented in order to help researchers build intuition and avoid pitfalls when moving from the context of graphs to the context of periodic graphs. One such misconception is that results on $c(G)$ carry over to $c(\tg)$ and we show this is not true for a simple extension of Berarducci and Intrigila's result in \cref{thm:cnt_retracts_copnum}.

In passing, we also prove general results. We start with results on \emph{temporal corners} (\cref{obs:temp_corner,prop:k_temp_corner}) whose contrapositives are used to build periodic graphs with high \copnum{s}, like in \cref{subsec:non_exist_temp_corners}. Then, we show an extension of Berarducci and Intrigila's bound that uses retracts (\cref{thm:retract_foremost_iso}). Finally, we prove a generalization of Joret et al's upper bound with the treewidth (\cref{thm:treewidth_foot_pg}). \Cref{thm:treewidth_foot_pg} presents the connection between separation properties of $G$ and $\ctgmax{G}$.

\section{Definitions}
\label{sec:definitions}

\subsection{Graphs and Time}

In what follows, we write $\Z^+$ for the set of positive integers including zero and $\Z_k$ for the set of integers modulo $k$. Given any integer $i$, we let $[i]_k$ be the integer in $\Z_k$ such that $i \equiv [i]_k \pmod{k}$.

\subsubsection{Static Graphs}
\label{sec:static}

We denote by $G = (V,E)$, or sometimes by $G = (V(G), E(G))$, the graph (directed or undirected) with set of vertices $V$ and set of edges $E$. We write $\duv{u}{v}$ for a \emph{directed} edge from $u$ to $v$ and $\uv{u}{v}$ for an \emph{undirected} one. A self-loop is an edge of the form $\duv{u}{u}$ or $\uv{u}{u}$. We say that $G$ is reflexive if every node has a self-loop. Unless stated otherwise, we consider all graphs to be undirected and reflexive. We say a graph $G'$ is a \emph{subgraph} of $G$, written $G'\subseteq G$, if $V(G')\subseteq V(G)$ and $E(G')\subseteq E(G)$. For any subset $V'\subseteq V$, we write $G[V'] = (V',E')$ for the subgraph of $G$ such that $E'$ contains all edges of $E$ that have both endpoints in $V'$ and say $G[V']$ is the subgraph of $G$ \emph{induced} by $V'$. We also write $G\setminus V'$ for the subgraph $G[V\setminus V']$ of $G$ and $G\setminus H$ for $G[V\setminus V(H)]$ when $H \subseteq G$. When $G$ is undirected, we write $N_G(u) := \set{v\in V; \uv{v}{u}\in E}\setminus\set{u}$ and $N_G[u] := N_G(u) \cup \set{u}$ for any node $u$. The \emph{degree} of $u$ in $G$, denoted $\deg_G(u)$, is given by $\abs{N_G(u)}$. Note that self-loops are not counted in $\deg_G(u)$.

For reasons apparent later, we shall refer to a graph $G$ so defined as a  {\em static graph}. 

\subsubsection{Temporal  Graphs}
 
A {\em time-varying graph} $\tg$ is a graph whose set of edges changes in time\footnotemark. A {\em temporal graph} is a time-varying graph where the set of time instants is $\Z^+$.
\footnotetext{
  The terminology in this section is mainly from \cite{Casteigts2013}.
}

A temporal graph $\tg$ is represented as an infinite sequence $\tg = (G_0, G_1, \dots)$ of static \emph{reflexive} and \emph{undirected} graphs $G_i =(V,E_i)$ on the same set of vertices $V$. We shall denote by $n=\abs{V}$ the number of  vertices of $\tg$. The graph $G_i$ is called the \emph{snapshot} of $\tg$ \emph{at time} $i\in \Z^+$ and the aggregate (undirected) graph $G = \left(V, \bigcup_{i} E_i\right)$ is called the {\em footprint} of $\tg$. We always write $V$ for the set of vertices of a temporal graph $\tg$.

Given two nodes $x, y\in V$, a \emph{journey}, from $x$ to $y$ starting at time $t$ is any finite sequence $\pi(x,y) = ((z_0, z_1),(z_1,z_2),\dots,(z_{k-1},z_k))$ where $z_0=x, z_k=y$, and  $\uv{z_i}{z_{i+1}}\in E(G_{t+i})$ for $0 \leq i < k$.  

A temporal graph $\tg$ is \emph{temporally connected} if for any $u,v\in V$ and any time $t\in \Z^+$ there is a journey from $u$ to $v$ that starts at time $t$. Observe that if $\tg$ is temporally connected, then its footprint is connected even when all its snapshots are disconnected. A temporal graph  $\tg$ is said to be \emph{always connected} (or $1$-{\em interval connected}) if all its snapshots are connected. 

A temporal graph $\tg$ is {\em periodic} if there exists a positive integer $p$ such that for all $i\in \Z^+$, $G_i = G_{i+p}$. If $p$ is the smallest such integer, then $p$ is called the {\em period} of $\tg$. We shall represent a periodic graph $\tg$ with period $p$ as $\tg = \peri{G_0,\dots, G_{p-1}}$. An example of a periodic graph $\tg$ with period $p=4$ is shown in \cref{fig:example_tg}. Observe that $\tg$ is temporally connected, however most of its snapshots are disconnected graphs. Note also that when $\tg$ is periodic, then it is temporally connected if and only if its footprint $G$ is connected. In this work we assume all periodic graphs are temporally connected unless specified otherwise.

Given a node $u\in V$ and a time $t\in \Z_p$, we write $\neig{t}{u}{\tg} := N_{G_t}[u]$ for the neighbourhood of $u$ at time $t$ and $\deg_t(u) := \abs{N_{G_t}(u)}$ for the degree of $u$ at time $t$. Furthermore, given a subset $V'\subseteq V$, the induced periodic subgraph $\tg[V']$ of $\tg$, induced by $V'$, is the periodic graph $\tg[V'] := \peri{G_0[V'], \dots, G_{p-1}[V']}$. If $H\subseteq G$, we write $\tg[H] := \peri{G_0[V(H)], \dots, G_{p-1}[V(H)]}$.

Let us point out the  obvious but useful fact that static graphs are periodic graphs with period $p=1$. 

\subsubsection{Arena}

All graphs in this paper are undirected, except for the following class of directed graphs we shall call \emph{arenas}. The notion of arena is similar to the notion of \emph{static expansion} in \cite{fluschnik2020time}.

\begin{definition}[Arena]\label{def:corners_arena}
Let $k\geq 1$ be any integer and $W$ be a non-empty finite set. An {\em arena} of length $k$ on $W$ is any static directed graph $\dg = (\Z_k\times W, E(\dg))$ where
$E(\dg) \subseteq \set{((i,w),([i+1]_{k}, w')); i\in \Z_k\mbox{ and } w,w'\in W}$.
\end{definition}  

\begin{figure}[h!]
\centering
\includegraphics[scale=1]{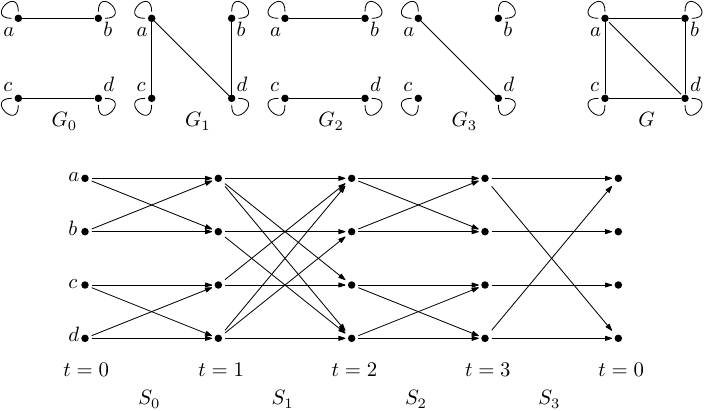}
\caption{A periodic graph $\tg = \peri{G_0,G_1,G_2}$ with its footprint $G$ and corresponding arena}
\label{fig:example_tg}
\end{figure}

A  periodic graph $\tg = \peri{G_0,\dots, G_{p-1}}$ with period $p$ and set of nodes $V$ has a unique correspondence with the  arena $\dg = (\Z_p\times V, E(\dg))$ where, for all $i\in \Z_p$, $((i,u),([i+1]_p,v))\in E(\dg) \iff \uv{u}{v}\in E(G_i)$, called the {\em arena of} $\tg$. The arena $\dg$ of $\tg$ explicitly preserves the snapshot structure of $\tg$ since there are obvious bijections between the snapshots of $\tg$ and subgraphs of $\dg$. An example of a periodic graph $\tg$ and its arena $\dg$ is shown in \Cref{fig:example_tg}. In the following, when no ambiguity arises, $\dg$ shall indicate the
arena of $\tg$.

The vertices of an arena $\dg$ will be called {\em temporal nodes}. Given a temporal node $(i,u)\in V(\dg)$ we shall denote by $\neig{i}{u}{\dg}$ the set of its outneighbours, and by $\npg{i}{u}{\dg}=\set{v\in V; ([i+1]_p,v)\in \neig{i}{u}{\dg}}$ the corresponding set of nodes in $G_{i}$. Given an arena $\dg = (\Z_p\times V, E)$ and a node $u\in V$, we write $\dg\setminus \set{u}$ for the arena $\dg' = (\Z_p \times V\setminus \set{u}, \set{\duv{(t,x)}{(t+1,y)}\in E; x\neq u\neq y})$.
 
\subsection{Cops \& Robber Game in Periodic Graphs}

\subsubsection{Basics}

The extension of the game of Cops and Robber from static to temporal graphs is quite natural. Initially, first the cops, then the robber, choose a starting position on the vertices of $G_0$. Then, at each time $t\in \Z^+$, first the cops, then the robber, move to vertices adjacent to their current positions in $G_{[t]_p}$. Thus, in round $t$, the players are in $G_{[t]_p}$ and, after making their moves, they find themselves in $G_{[t+1]_p}$ in the next round. The game ends and the cops win if and only if at least one cop moves to the vertex currently occupied by the robber. The robber wins by forever preventing the cops from winning. 

\subsubsection{Configurations and Strategies}

Let $k\geq 1$ cops play on $\tg$. A configuration is a pair of possible positions for $k$ cops and the robber when the game is played on $\tg$, written as $\conf((t,c_1,\dots,c_k),(t',r))$ where $t$ and $t'\in\set{t-1,t}$ are times, $c_1,\dots,c_k$ are the positions of the cops and $r$ is the position of the robber. A strategy for the cops is a function $\cs$ that maps each configuration to a new position for the cops and robber strategies are similarly defined. We say a cops strategy $\cs$ is \emph{feasible} if whenever $\cs((t,c_1,\dots,c_k),(t,r)) = (t+1,c_1',\dots,c_k')$, then $c_i'\in \neig{t}{c_i}{\tg}$ for every $1\leq i\leq k$. The same holds for robber strategies. A configuration $\conf((t,c_1,\dots,c_k),(t,r))$ is said to be $k$-\emph{copwin} if there exists a strategy $\cs$ such that, starting from $\conf((t,c_1,\dots,c_k),(t,r))$, the cops win the game regardless of the strategy used by the robber. Moreover, $\cs$ is said to be $k$-copwin on $\tg$ if there exists $k$ nodes $u_1,\dots,u_k$ such that for any node $v$, $\cs$ is winning from $\conf((0,u_1,\dots,u_k),(0,v))$.

We say $\tg$ and $\dg$ are $k$-copwin if the cops have a $k$-copwin strategy on $\tg$. When $k = 1$, we write copwin instead of $1$-copwin.

We say a cops strategy on $\tg$ is \emph{stubborn} if it describes where the cops should go on $G$ and the cops make their moves as their incident edges become available in $\tg$. When the cops follow such a strategy, we say they move \emph{stubbornly}.

\subsubsection{Temporal Corners}

\begin{figure}
\centering
\includegraphics[scale = 1]{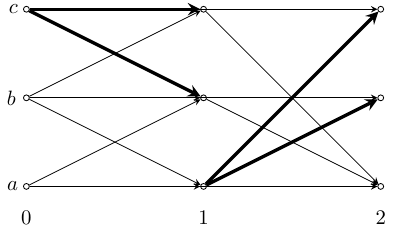}
\caption{An arena with a temporal corner $(0,c)$ of $(1,a)$}
\label{fig:temp_corner_exa}
\end{figure}

A temporal node $(t,u)$ in an arena $\dg$ is a \emph{temporal corner} of a \emph{temporal cover} $(t+1,v)$ if $u\neq v$ and
\[
\npg{t}{u}{\dg}\subseteq \npg{t+1}{v}{\dg}.
\]
In $\dg$, every time it moves, a single cop ends its turn in the snapshot ahead of the robber. Thus, this definition of temporal corner encapsulates the usual meaning of corner that \emph{after the cop has moved, no matter where the robber plays, the robber gets captured the next time the cop moves.} We add the restriction $u\neq v$ since, if the robber stands at $(t,u)$, when the cop moves to $(t+1,u)$, the game ends with the cop winning. This situation is different from the situation when the cop moves to a temporal cover of the robber's position: then the robber is still not captured, but cannot escape from the cop. The simplest way to determine if an arena $\dg$ contains a temporal corner is to iterate through all temporal nodes $(t,u)$ and nodes $v$ and verify whether the inclusion $\npg{t}{u}{\dg}\subseteq \npg{t+1}{v}{\dg}$ holds or not. However, since $u\in \npg{t}{u}{\dg}$, we can limit our search to those nodes $v$ such that $u\in\npg{t+1}{v}{\dg}$. As an example, on the arena $\dg$ of \cref{fig:temp_corner_exa}, $(0,c)$ is a temporal corner of $(1,a)$ because $\set{b,c} = \npg{0}{c}{\dg} \subseteq \npg{1}{a}{\dg} = \set{a,b,c}$. The relevant edges are drawn in bold.

More generally, we say $(t,x)$ is a \emph{$k$-temporal corner} of $(t+1,y_1),\dots,(t+1,y_k)$, if $x\notin \set{y_1,\dots,y_k}$ and
\[
\npg{t}{x}{\dg} \subseteq \bigcup_{i=1}^k \npg{t+1}{y_i}{\dg}.
\]

The following two results relate $k$-temporal corners to $k$-copwin periodic graphs.

\begin{lemma}\label{obs:temp_corner}
Every copwin arena contains a temporal corner.
\end{lemma}
\begin{proof}
Let $\dg$ be a copwin arena. Observe that if the cop wins in a single move in $G_0$, then $G_0$ has a universal vertex $u$ and every temporal node $(p-1,x)$ is a temporal corner of $(0,u)$. Otherwise, consider a configuration $\conf((t+1,v),(t,u))$ such that wherever the robber moves to, the robber gets captured by the cop. This configuration exists because $\dg$ is copwin. Since it is the robber's turn to play, for every $w\in \npg{t}{u}{\dg}$, there exists a $z\in \npg{t+1}{v}{\dg}$ such that $z = w$. In other words, $\npg{t}{u}{\dg}\subseteq \npg{t+1}{v}{\dg}$ and $(t,u)$ is a temporal corner of $(t+1,v)$. 
\qed
\end{proof}

\begin{proposition}\label{prop:k_temp_corner}
Every $k$-copwin arena contains a $k$-temporal corner.
\end{proposition}
\begin{proof}
Let $\dg$ be a $k$-copwin arena. If the cops win in a single move in $G_0$, then $G_0$ has a dominating set $\set{v_1,\dots,v_k}$ of size $k$ and every temporal node $(p-1,x)$ is a $k$-temporal corner of $(0,v_1),\dots,(0,v_k)$. Otherwise, consider a configuration $\conf((t+1,y_1,\dots,y_k),(t,x))$ such that wherever the robber moves to, the robber gets captured by the cops. Since it is the robber's turn to play, for every $w\in \npg{t}{x}{\dg}$, there exists a node $y_i$ and a neighbour $z\in \npg{t+1}{y_i}{\dg}$ such that $z = w$. In other words, $\npg{t}{x}{\dg} \subseteq \bigcup_{i=1}^k \npg{t+1}{y_i}{\dg}$. Therefore, $\dg$ contains a $k$-temporal corner.
\qed
\end{proof}

\Cref{prop:k_temp_corner} implies that if $\tg$ does not contain any $k$-temporal corner then it cannot be $k$-copwin. This contrapositive will often come in handy in the next sections. However, the converse is not true in general.

\section{Comparing the cop number of a periodic graph with the cop numbers of its footprint}

\subsection{A motivational example}

We define $\ctgmax{G}$ as the maximal \copnum{} $c(\tg)$ of any periodic graph $\tg$ with footprint $G$. 

Recall that a periodic graph $\tg$ is defined as a sequence of (possibly disconnected) subgraphs of a graph $G$ on the same set of nodes $V$. Because those subgraphs define the structure of $\tg$, it seems natural to wonder if $\ctgmax{G}\leq f(c(G))$ for some function $f$. We start with the simplest function and inquire if $\ctgmax{G}\leq c(G)$. If so, then we would have a \say{simple} upper bound on $\ctgmax{G}$. Computing $c(G)$ is \expc{} in general (see \cite{Kinnersley2015}). Thus, not only would this bound still be hard to compute if it held, but if it were not true, then $\ctgmax{G}$ might be much more difficult to compute in general. 

Unfortunately, we answer this question in the negative. The following construction is based on the hypercube $Q_3$. It is known that $c(Q_k) = \cl{\frac{k+1}{2}}$ \citep{maamoun1987game}. We use the usual construction for $Q_3$ and write each node as a bit sequence of length $3$.
\begin{figure}[h!]
\centering
\includegraphics{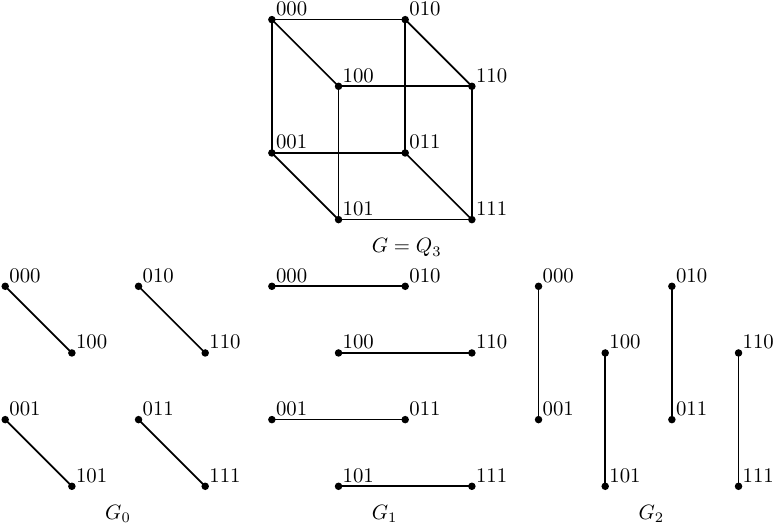}
\caption{Periodic graph used in \cref{prop:hypercubeQ3_3copwin}
}
\label{fig:hypercubeQ3_3copwin}
\end{figure}

\begin{proposition}[3-copwin cube periodic graph]\label{prop:hypercubeQ3_3copwin}
There exists a periodic graph $\tg$ with footprint $G = Q_3$ such that $c(\tg) = 3$. 
\end{proposition}
\begin{proof}
Let us describe a periodic graph $\tg = \peri{G_0,G_1,G_2}$ with footprint $G = Q_3$ such that $c(\tg) = 3$. In $G_0$, only edges that change the first bit appear. In $G_1$, only those that change the second bit appear and so on for $G_2$. This is shown in \cref{fig:hypercubeQ3_3copwin}. We claim two cops cannot catch the robber. 

Let us call a $4$-cycle of $Q_3$ a \emph{face}. Any subgraph of any snapshot of $\tg$ that induces a face in $Q_3$ is also called a face.   

We wish to preserve the following invariant for the robber:
\begin{itemize}
\item[(I)]
Two cops are never in the same face as the robber at the beginning of their turn.
\end{itemize}

Let us show that if invariant $I$ holds in $G_t$ for all $t\geq 0$, then no matter where the cops move, the robber will not be on a $2$-temporal corner of the cops' positions. Then, by \cref{prop:k_temp_corner}, $c(\tg)>2$. 

Suppose otherwise, that $(t,u)$ is a $2$-temporal corner of $(t+1,v), (t+1,w)$. Observe that for any time $t$ and nodes $r,c\in V$ with $r\neq c$, $\abs{\neig{t}{r}{\tg} \cap \neig{t+1}{c}{\tg}}\leq 1$. Let $\neig{t}{u}{\tg} = \set{x,u}$. Moreover, $\deg_t(u) = \deg_{t+1}(v) = \deg_{t+1}(w) = 1$ and $E(G_t)\cap E(G_{t+1}) = \emptyset$. Thus, either $\neig{t+1}{v}{\tg} = \set{v,x}$ and $\neig{t+1}{w}{\tg} = \set{w,u}$, or $v = x$ and $\neig{t+1}{w}{\tg} = \set{w,u}$, without loss of generality. In both cases, $(uxvw)$ is (part of) a face of $G$. Indeed, only the edges that change the same bit appear at time $t+1$, so either $\uv{w}{v}\in E(G)$ (when $v\neq x$) or $\uv{w}{y}\in E(G)$ where $\neig{t+1}{v}{\tg} = \set{v,y}$ (when $v = x$). Then, we either have $\uv{x}{u}, \uv{x}{v}, \uv{v}{w}, \uv{w}{u}\in E(G)$ or $\uv{x}{u}, \uv{v}{y}, \uv{y}{w}, \uv{w}{u}\in E(G)$ which are both faces of $Q_3$.

Therefore, no matter where the robber moves to, the cops will start their turn in $G_{t+1}$ in the same face as the robber. This contradicts our invariant $I$.

Let us show now that the robber can play so that $I$ is always true. 

Since $Q_3$ has $6$ faces, the robber can avoid choosing an initial position in a face that contains both cops in $G_0$.

Let us show that if $I$ is true before the cops have played, then it will remain true after the robber has moved. Without loss of generality because of the symmetries in the snapshots of $\tg$, suppose the robber is on $(0,000)$ and the cops are on $c_1\neq 000\neq c_2$ in $G_1$. We moreover assume that $\set{000,c_1,c_2}$ is contained in a face of $Q_3$ since otherwise the robber could easily move so that $I$ is true. 

The node $000$ is contained in $3$ faces of $Q_3$, but no cop is on $000$, so there are $9$ cases.
\begin{enumerate}
\item 
The cops are on $\uv{(001)}{(011)}$. In $G_0$, the robber moves to $100$ and avoids ending in the same face as the cops.
\item 
The cops are on $\uv{(001)}{(101)}$. The only way for the cops to occupy this edge in $G_1$ was for them to start on this edge in $G_0$ and stay on their position. Therefore, in $G_0$ the robber is in the same face $(000,100,101,001)$ as the cops before the cops have played. This violates our assumption that $I$ was true before the cops played.
\item 
The cops are on $\uv{(010)}{(011)}$. In $G_0$, the robber moves to $100$ and avoids ending in the same face as the cops.
\item 
The cops are on $\uv{(010)}{(110)}$. The only way for the cops to occupy this edge in $G_1$ was for them to start on this edge in $G_0$ and stay on their position. Therefore, in $G_0$ the robber is the same face $(000,100,110,010)$ as the cops. This violates our assumption that $I$ was true before the cops played.
\item 
The cops are on $\uv{(100)}{(101)}$. Before the cops moved in $G_0$, all players were in the same face $(000,100,101,001)$. This violates our assumption that $I$ was true before the cops played.
\item 
The cops are on $\uv{(100)}{(110)}$. Before the cops moved in $G_0$, all players were in the same face $(000,100,110,010)$. This violates our assumption that $I$ was true before the cops played.
\item
The cops are on $\set{010, 100}$. Before the cops moved in $G_0$, all players were in the same face $(000,010,110,100)$. This violates our assumption that $I$ was true before the cops played.
\item 
The cops are on $\set{010, 001}$. The robber moves to $100$ and avoids ending in the same face as the cops.
\item
The cops are on $\set{100, 001}$. Before the cops moved in $G_0$, all players were in the same face $(000,100,101,001)$. This violates our assumption that $I$ was true before the cops played.
\end{enumerate}
Thus, from $(0,000)$ the robber can move so that $I$ will be true in $G_1$. By symmetry, this holds for every robber position under optimal play. Therefore, the invariant $I$ always holds and $c(\tg)>2$. 

We argue that $c(\tg)\leq 3$. Let three cops start on $000$, $010$ and $111$ in $G_0$. The robber starts either on $001$ or $101$ in order to avoid getting captured in the first turn. The cops stay still in $G_0$. The robber must end its turn on $001$, otherwise the cop on $111$ would make the catch in $G_1$. In $G_1$, no cops move so that the robber end its turn either on $001$ or $011$. Since there are cops on $000$ and $010$, the robber gets captured in $G_2$.
\qed \end{proof}

Notice that in the previous result, $c(\tg) = 3 > c(G) = 2 = \gamma(G)$, where $\gamma(G)$ is the domination number of $G$. In the static case, the domination number is a trivial upper bound on $c(G)$.

In fact, this also shows another interesting result. Baird et al. \cite{baird2013minimum} confirmed that in the static case, the smallest $3$-copwin graph is the Petersen graph with $10$ vertices and Turcotte and Yvon (\cite{turcotte20214}) later pushed the investigation further by showing that $4$-copwin graphs have at least $19$ vertices. \Cref{prop:hypercubeQ3_3copwin} shows that the smallest $3$-copwin \emph{periodic} graph has at most $8$ vertices. One can show it cannot have four vertices or less. This shows that if one fixes a particular \copnum{} $c$, it is possible to lower the number of nodes required to get a periodic graph with \copnum{} $c$ compared to the static case. 
\begin{theorem}
The smallest $3$-copwin periodic graph has at least five vertices and at most eight.
\end{theorem}

\subsection{Retracts and Tree Decompositions}

\begin{figure}
\centering
\begin{subfigure}[t]{0.4\linewidth}
\includegraphics[scale=0.9]{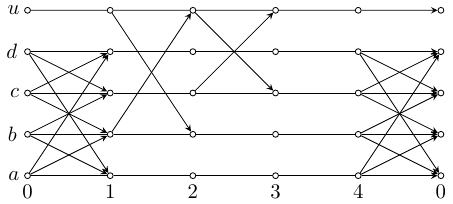}
\caption{Arena $\dg$}
\label{subfig:retract_cnt}
\end{subfigure}
\hfill
\begin{subfigure}[t]{0.4\linewidth}
\includegraphics[scale=0.9]{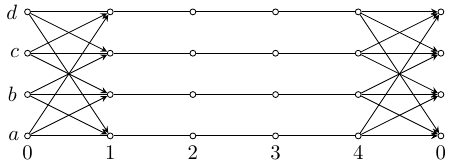}
\caption{Arena $\dg\setminus \set{u}$}
\label{subfig:retract_cnt_minusu}
\end{subfigure}
\caption{$c(\tg[\set{a,b,c,d}]) > c(\tg)$}
\label{fig:retract_cnt}
\end{figure}

An important tool in the study of Cops and Robber games is the concept of retract. A \emph{retraction} $h$ of a graph $G$ is a homomorphism from $G$ to one of its subgraphs $H$ that is the identity on $H$. The subgraph $H$ is called a \emph{retract} of $G$. Retracts are used to show that copwin graphs are dismantlable (see for example \cite{Bonato2011f}), find a tighter upper bound on $c(G)$ with block decompositions \cite{hill2008cops} and even to help showing that the Petersen graph is the smallest $3$-copwin graph \cite{baird2013minimum}. One important result on retracts is a theorem of \citet{berarducci1993cop} stating that $c(H)\leq c(G)$ when $H$ is a retract of $G$. In other words, upper bounds on $c(G)$ are carried over when taking retracts.
Unfortunately, this does not hold true in general for periodic graphs when considering retracts of the footprint.
\begin{proposition}\label{thm:cnt_retracts_copnum}
There exists a periodic graph $\tg$ with footprint $G$ and a retract $H$ of $G$ such that $c(\tg[H]) > c(\tg)$.
\end{proposition}
\begin{proof}
Let $G$ be the graph with vertex set $\set{a,b,c,d,u}$, $4$-cycle $(a,b,c,d)$ and $3$-clique $(b,u,c)$. Then, $H = C_4$ is a retract of $G$ obtained by mapping $u$ to either $b$ or $c$. Let us construct a periodic graph $\tg$ on $G$ such that $c(\tg) = 1$ while $c(\tg[H]) = 2$.

Consider the arena $\dg$, shown in \cref{subfig:retract_cnt}, that corresponds to $\tg$. If the cop starts on $(0,a)$ it can capture the robber no matter where the latter starts. Thus, this is copwin. The arena $\dg\setminus \set{u}$, that corresponds to $\tg[H]$ and shown in \cref{subfig:retract_cnt_minusu}, contains two identical $4$-cycles, so it cannot be copwin. However, two cops can win on $\dg\setminus \set{u}$. Therefore, $c(\tg) = 1 < c(\tg[H]) = 2$.
\qed \end{proof}

If we add one more assumption, we can recover Berarducci and Intrigila's original result.

\begin{theorem}\label{thm:retract_foremost_iso}
Let $\tg$ be a periodic graph with footprint $G$ and $h : G\to H$ be a retraction of $G$. If for every time $t$, $h(G_t)$ is a retract of $G_t$, then
\[
c(\tg[H])\leq c(\tg).
\]
\end{theorem}
\begin{proof}
For every time $t$, we write $H_t = h(G_t)$. For any edge $\uv{x}{y}\in E(G_t)$, $\uv{h(x)}{h(y)}\in E(H_t)$ because $h$ is a retraction on the snapshots. 

Let $\cs$ be any strategy for $k$ cops on $\tg$. Let $C = (c_1,\dots,c_k)$ be the sequence of positions occupied by the cops at any time $t$. If $\cs((t,C), (t,r)) = (t+1,C')$ in $\tg$, with $C' = (c_1',\dots,c_k')$, then $\uv{c_i}{c_i'}\in E(G_t)$ for every $1\leq i\leq k$. By the above argument, $\uv{h(c_i)}{h(c_i')}\in E(H_t)$. Thus, every cops strategy $\cs$ on $\tg$ has a corresponding strategy $\cs^h$ on $\tg[H]$ such that if $\cs((t,C),(t,r)) = (t+1,C')$, then $\cs^h$ moves the cop on $c_i$ from $h(c_i)$ to $h(c_i')$ for every $1\leq i\leq k$.

Thus, let $k = c(\tg)$ and let the cops play a $k$-winning strategy $\cs$ on $\tg$ while the robber is restricted to play on $\tg[H]$. Since $\tg$ is $k$-copwin, the cops eventually move to a position $(t+1,c_1),\dots,(t+1,c_k)$ while the robber is on $(t,r)$ such that no matter where the robber moves to it gets captured on the next turn. That is, for every $x\in \neig{t}{r}{\tg[H]}$, there exists $\uv{x}{c_i}\in E(G_{t+1})$. Then, by the above argument, $\uv{h(x)}{h(c_i)}\in E(H_{t+1})$. Thus, $\conf((t+1,h(c_1), \dots, h(c_k)), (t,r))$ is a winning configuration for the cops in $\tg[H]$. As we argued, this configuration can be reached when the cops follow the strategy $\cs^h$ that corresponds to $\cs$ in $\tg[H]$. Therefore, the robber gets captured by the $k$ cops in $\tg[H]$.

Thus, $c(\tg[H])\leq c(\tg)$.
\qed \end{proof}
Consider again the arena in \cref{fig:retract_cnt}. Let us show why \cref{thm:retract_foremost_iso} does not apply in this case. Let $h : G\to H$ be any retraction. The node $u$ can only be mapped to either $b$ or $c$ because $N_G(u) = \set{b,c}$. Suppose that $h(u) = b$. Then, $\uv{u}{c}\in E(G_2)$ implies that $\uv{h(u)}{h(c)} = \uv{b}{c}\in E(H_2)$, which is not the case. Similarly, if $h(u) = c$, then $\uv{u}{b}\in E(G_1)$ implies that $\uv{h(u)}{h(b)} = \uv{c}{b}\in E(H_1)$, which is also not the case. Thus, $h$ cannot be a retraction of every snapshot.

Before presenting the statement of \cref{thm:treewidth_foot_pg} (see below), we need to define the \emph{tree decomposition} and the \emph{treewidth} of a graph. We refer the reader to \citet{Diestel2000} for more details.
\begin{definition}[Tree decomposition]
A \emph{tree decomposition} of a graph $G$ is a pair $(T, \tb=\set{B_x; x\in V(T)})$ where $T$ is a tree, $\tb$ is a family of subsets of $V(G)$ called bags and such that:
\begin{itemize}
\item 
$\bigcup_{x\in V(T)} B_x = V(G)$;
\item
For every edge $uv\in E(G)$, there exists some $x\in V(T)$ such that $u,v\in B_x$;
\item
For every vertex $u\in V(G)$, the set $\set{x\in V(T); u\in B_x}$ induces a subtree of $T$.
\end{itemize}
The \emph{width} of a tree decomposition $(T, \tb)$ is $\max_{x\in V(T)} \abs{B_x} - 1$. The \emph{treewidth} of $G$, written $\tw(G)$, is the minimum width among all tree decompositions of $G$. We say a tree decomposition is \emph{minimal} if its width equals $\tw(G)$.
\end{definition}

A tree decomposition $(T, \tb)$ of $G$ with width $k$ is \emph{smooth} if every bag $B_x\in \tb$ has size $k + 1$ and for every edge $xy\in E(T)$, $\abs{B_x\cap B_y} = k$. Any tree decomposition of $G$ can be transformed into a smooth tree decomposition with the same width (see \cite{bodlaender1996linear}). Smooth tree decompositions are similar to \emph{normalized} tree decompositions (\cite{harvey2017parameters}).

Joret et al. \cite{joret2008cops} proved that $c(G)\leq \tw(G)/2 + 1$. We prove the following.
\begin{theorem}\label{thm:treewidth_foot_pg}
For every graph $G$, $\ctgmax{G} \leq \tw(G)+1$.
\end{theorem}
\begin{proof}
Let $G$ have treewidth $k$ and $(T, \tb)$ be a minimal tree decomposition of $G$ that is \emph{smooth}. Recall that for any $\uv{x}{y}\in E(T)$, $B_x\cap B_y$ is a cutset of $G$. Also, we write $T_{x,y}$ for the subtree of $T\setminus \set{x}$ that contains $y$.

Let $\tg = \peri{G_0,\dots, G_{p-1}}$ be any periodic graph with footprint $G$. For any bag $B_x$ of $T$, $k+1$ cops can be positioned in $B_x$ so that if the robber moves into $B_x$ at any time it gets captured immediately. Thus, let $k+1$ cops start on the $k+1$ nodes of any bag $B_{x_0}$ in $G_0$. Let the robber start on some node $r_0$. 

Suppose the robber has not been captured and let $x_1$ be the unique neighbour of $x_0$ in $T$ such that $r_0$ is in a bag of $T_{x_0,x_1}$. Because $T$ is smooth, $\abs{B_{x_0}\cap B_{x_1}} = k$ and there exists a unique node $v\in B_{x_1} \setminus B_{x_0}$ and a unique node $u\in B_{x_0} \setminus B_{x_1}$. 

Recall that $\tg$ is temporally connected, so let the cop on $u$ walk to $v$, traversing edges whenever possible.

Meanwhile, all other cops stay still. In order to escape from $T_{x_0,x_1}$, the robber has to move trough a node of $B_{x_0}\cap B_{x_1}$ because this is a cutset of $G$. But, the nodes in this set will all be occupied while the travelling cop moves to $v$. Therefore, once this cop arrives on $v$, the robber is still in $T_{x_0,x_1}$. Furthermore, at that time all nodes of $B_{x_1}$ are occupied by a cop.

When this happens, the robber's territory is reduced. It follows by successively applying this reasoning that $c(\tg)\leq k + 1$.
\qed \end{proof}

\section{About the maximum and minimum cop numbers of the snapshots}

In the previous section, we showed an example where $c(\tg)>c(G)$. Here, we take the opposite direction and show examples where $c(\tg)<c(G)$, to highlight that both cases are possible. The main question here is: How low can $c(\tg)$ be compared to $c(G)$? 

The following result, along with its corollary, answers this.
\begin{lemma}\label{lem:kkpkp_copwin}
For any $k\geq 3$ and $1\leq k'<k$ such that there exists a $k$-copwin graph with a spanning $k'$-copwin subgraph, there exists an at most $k'$-copwin periodic graph whose footprint is $k$-copwin.
\end{lemma}
\begin{proof}
Let $G$ be any graph with $c(G) = k$ and $H$ a spanning $k'$-copwin subgraph of $G$. It suffices to let $H$ appear long enough in a periodic graph $\tg$ so that $c(\tg) \leq k'$. Then, we can cover the remaining edges of $G$ in the remaining snapshots with spanning trees, so that $c(G) = k$ and $c(\tg)\leq k'$. 
\qed \end{proof}

Let $H = (V(H), E(H))$ be a graph and $d_H(u,v)$ be the distance in $H$ from $u\in V(H)$ to $v\in V(H)$. We define $\rad(H):= \min_{x\in V(H)}\max_{y\in V(H)} d_H(x,y)$.

\begin{corollary}\label{cor:311_pg}
There exists a copwin periodic graph whose footprint is $3$-copwin.
\end{corollary}
\begin{proof}
Let $G$ be the Petersen graph, so $c(G) = 3$. Let $\tg = \peri{G_0,\dots, G_{p-1}}$ be such that the first $G_0,\dots,G_l$ snapshots contain a minimum spanning tree $T$ of $G$, with $l\geq \rad(T)$. The remaining snapshots contain different spanning trees to cover the edges of $G$. Therefore, $c(\tg) = 1$. 
\qed \end{proof}

From \cref{lem:kkpkp_copwin}, it follows that $c(\tg)$ can be arbitrarily small compared to $c(G)$ since $\tg$ might be much sparser than its footprint.
Then, one can wonder what happens if we limit the possible \copnum{s} of all snapshots. For this purpose, given a periodic graph $\tg = \peri{G_0,\dots, G_{p-1}}$ we define $\Sdo:=\argmin_{0\leq i\leq p-1} c(G_i)$ as well as $\Sup:=\argmax_{0\leq i\leq p-1}c(G_i)$. In the next results, we study whether $c(G), c(\Sdo)$ and $c(\Sup)$ relate to $c(\tg)$ and, if so, how.

\begin{figure}
\centering
\foreach \i in {0,...,8}{%
    \begin{subfigure}{0.3\textwidth}
    \includegraphics[scale=0.5]{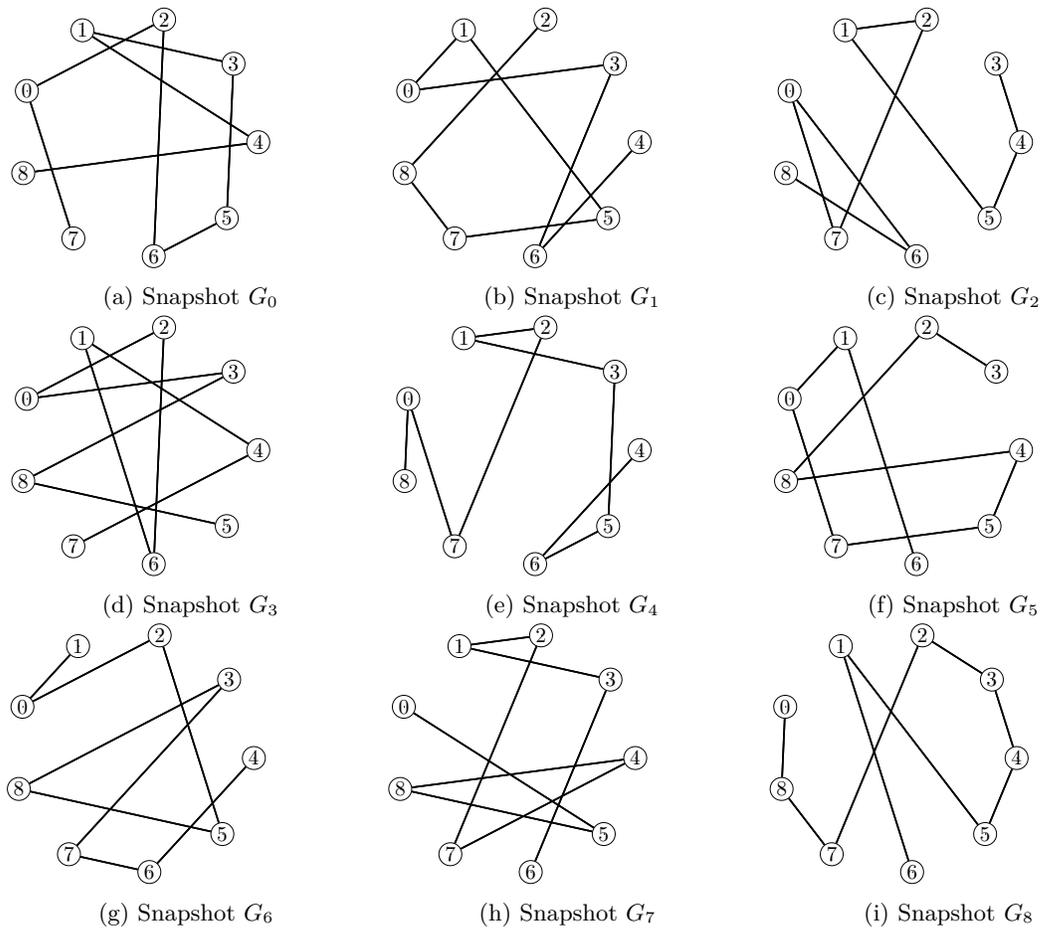}
    \hfill
    \protect\caption{Snapshot $G_{\i}$}
    \end{subfigure}
}
\caption{The $(1,1,2)$-copwin periodic graph presented in \cref{conj:112_periodicgraph}}
\label{fig:pg_thm112}
\end{figure}

\begin{figure}
\centering
\includegraphics[scale=0.5]{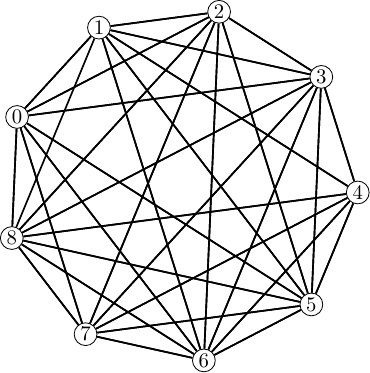}
\caption{The footprint of the periodic graph from \cref{conj:112_periodicgraph}}
\label{fig:pg_thm112_foot}
\end{figure}

One might expect the following pair of inequalities to hold:
$
\min(c(G), c(\Sdo))
\leq
c(\tg) 
\leq 
\max(c(G), c(\Sup)).
$
Indeed, this forms the widest range of \copnum{s} that only uses parameters from $c(G_0), c(G_1),\dots,c(G_{p-1})$ along with $c(G)$. Nevertheless, we show both are false in \cref{conj:112_periodicgraph} and \cref{thm:up_low_bounds_impos}. These two results present examples that will have a larger importance in \cref{sec:compl_table}.
In the statement of \cref{conj:112_periodicgraph} and \cref{thm:up_low_bounds_impos}, we say that \emph{$\tg$ is $(a,b,c)$-copwin} if $c(G) = a, c(\Sup) = b$ and $c(\tg) = c$. This notation will be extensively used in the next section.
\begin{theorem}[(1,1,2)-copwin]\label{conj:112_periodicgraph}
The inequality $c(\tg) \leq \max(c(G), c(\Sup))$ is false.
\end{theorem}
\begin{proof}
Consider the periodic graph $\tg$ whose snapshots are shown in \cref{fig:pg_thm112}. The footprint $G$ is shown in \cref{fig:pg_thm112_foot}. Each snapshot is a path on $9$ vertices, so it is copwin. The footprint $G$ has a universal vertex (node $8$), thus it is also copwin. Inspection shows that $\tg$ contains no temporal corner, so it cannot be copwin by \cref{obs:temp_corner}. 

We now show that $\tg$ is $2$-copwin by describing a winning strategy for two cops. Let two cops start on nodes $2$ and $4$ in $G_0$. The robber must start on either $3,5$ or $7$ to avoid capture in $G_0$. 
\begin{enumerate}
\item
If the robber starts on $3$, the cop on $4$ moves to $1$ and the cop on $2$ moves to $0$. The robber cannot move to $1$ nor to $5$, since $\uv{1}{5}\in E(G_1)$, so it stays on $3$. In $G_1$, $\uv{0}{3}$ appears and one cop makes the catch.
\item
If the robber starts on $5$, the cop on $4$ moves to $1$ and the cop on $2$ moves to $6$. The robber cannot move to $6$ nor to $3$ because $\uv{6}{3}\in E(G_1)$. Therefore, the robber stays on $5$ and gets captured by the cop on $1$ that moves along the edge $\uv{1}{5}$ in $G_1$. 
\item 
If the robber starts on $7$, the cops move to $0$ and $8$ in $G_0$. The robber cannot move to $0$, so it stays on $7$ and gets captured by the cop on $8$ in $G_1$.
\end{enumerate}
\qed \end{proof}

\begin{theorem}[(2,2,1)-copwin]\label{thm:up_low_bounds_impos}
The inequality $\min(c(G), c(\Sdo))\leq c(\tg)$ is false.
\end{theorem}
\begin{proof}
Let $G$ be a bow tie graph formed with two $4$-cycles joined on a vertex $v$. Let $G_0=G_1=G_{2}$ be the subgraph of $G$ induced by removing one edge of a cycle and $G_{3}=G_4=G_{5}$ be the subgraph of $G$ induced by removing one edge of the other cycle. Clearly, $G$, $G_1,\dots,G_5$ are all $2$-copwin. However, on each snapshot there is a path. Let $\tg=\peri{G_0,\dots, G_{5}}$. Place one cop on $v$ in $G_0$. The robber initially starts on the cycle of $G_0$. The cop waits until $G_{3}$ for the path to appear underneath the robber. The robber cannot have crossed onto the other cycle since the cop is guarding the cutvertex $v$. The robber is now on a path of length at most $3$, so the cop has enough time to capture the robber.
\qed \end{proof}
\Cref{thm:up_low_bounds_impos} is a dual result of \cref{conj:112_periodicgraph} and together they serve to highlight how loose the connection is between the \copnum{s} of a periodic graph, its snapshots and footprint. 

\section{Completing the table of copwin periodic graphs}\label{sec:compl_table}

In the previous section, we showed that simple lower and upper bounds on $c(\tg)$ that depend on $c(G), c(G_0),\dots,\\ c(G_{p-1})$ do not hold in general. In this section, we want to further emphasize the disconnect between those values by showing that nearly all combinations of values (between $1$ and $3$) of $c(\tg)$, $c(G)$ and $c(\Sup)$ are possible (refer to \cref{tab:summ_results}).

In this section, we show that out of the $27$ possible combinations of those parameters, $3$ remain to be determined. We suspect that those would have to be discovered by computer search if they exist. For example, the $(1,1,2)$-copwin periodic graph, \cref{conj:112_periodicgraph}, could be extended to generate a $(2,1,3)$ or $(3,1,3)$-copwin periodic graph.

\begin{table}[ht]
\caption{Summary of results on periodic graphs. These are existence results of periodic graphs $\tg$ with \copnum{} $c(\tg)$. The \copnum{} of the footprint is noted $c(G)$ and the maximum \copnum{} of the snapshots is $c(\Sup)$}
\newcolumntype{R}{>{\raggedleft\arraybackslash}X}%
\begin{tabularx}{\textwidth}{XXXl @{\hspace{64\tabcolsep}} XXXl}
\toprule
$c(G)$ & $c(\Sup)$ & $c(\tg)$ & Reference & $c(G)$ & $c(\Sup)$ & $c(\tg)$ & Reference\\
\midrule
$1$ & $1$ & $1$ & \cref{prop:leftover_results} & $2$ & $2$ & $3$ & \cref{prop:leftover_results}\\
$1$ & $1$ & $2$ & \cref{conj:112_periodicgraph} & $2$ & $3$ & $1$ & \cref{prop:231_copwin}\\
$1$ & $1$ & $3$ & \mbox{\emph{Undetermined}} & $2$ & $3$ & $2$ & \cref{prop:leftover_results}\\
$1$ & $2$ & $1$ & \cref{prop:leftover_results} & $2$ & $3$ & $3$ & \cref{prop:leftover_results}\\
$1$ & $2$ & $2$ & \cref{lem:copw_footprint_2sn_noncopwin} & $3$ & $1$ & $1$ & \cref{cor:311_pg}\\
$1$ & $2$ & $3$ & \cref{prop:123-copwin_generic} & $3$ & $1$ & $2$ & \cref{prop:leftover_results}\\
$1$ & $3$ & $1$ & \cref{prop:leftover_results} & $3$ & $1$ & $3$ & \mbox{\emph{Undetermined}}\\
$1$ & $3$ & $2$ & \cref{lem:132_copwin} & $3$ & $2$ & $1$ & \cref{lem:321_copwin}\\
$1$ & $3$ & $3$ & \cref{prop:leftover_results} & $3$ & $2$ & $2$ & \cref{prop:leftover_results}\\
$2$ & $1$ & $1$ & \cref{prop:leftover_results} & $3$ & $2$ & $3$ & \cref{prop:leftover_results}\\
$2$ & $1$ & $2$ & \cref{prop:leftover_results} & $3$ & $3$ & $1$ & \cref{prop:leftover_results}\\
$2$ & $1$ & $3$ & \mbox{\emph{Undetermined}} & $3$ & $3$ & $2$ & \cref{prop:leftover_results}\\
$2$ & $2$ & $1$ & \cref{thm:up_low_bounds_impos} & $3$ & $3$ & $3$ & \cref{prop:leftover_results}\\
$2$ & $2$ & $2$ & \cref{prop:leftover_results} & & & &\\
\bottomrule
\end{tabularx}
\label{tab:summ_results}
\end{table}

First, let us show that from any particular result, we can always increase the number of nodes. This way, our constructions will remain general. 
\begin{lemma}\label{prop:abccop_larger_abccop}
For any $(a,b,c)$-copwin periodic graph $\tg$, with $n$ nodes and period $p\geq 2$, and any integer $N\geq n$, there exists a periodic graph $\tg'$ with $N$ nodes and period $p$ that is $(a,b,c)$-copwin. 
\end{lemma}
\begin{proof}
Let $\tg = \peri{G_0, \dots, G_{p-1}}$ with footprint $G$ be as in the statement and $P = (u_1,\dots,u_{N - n + 1})$ be a path with $N - n + 1$ nodes, where $N\geq n$ is some integer. Let $u\in V(G)$. For every $1\leq i\leq p-1$, let $H_i$ be obtained by identifying $u$ with $u_1$ and let $\tg' = \peri{H_0, \dots, H_{p-1}}$ with footprint $H$. By construction, $c(\tg'[G]) = c(\tg)$.

The robber on $\tg'$ can play on $\tg'[G]$ and win against less than $c(\tg'[G]) = c(\tg)$ cops. Therefore, $c(\tg') \geq c(\tg'[G]) = c(\tg)$.

Let us show that $c(\tg') \leq c(\tg)$. By construction, the homomorphism $h : H\to G$ that maps $V(P)$ to $u$ and is the identity on $G$ is a retraction of $H$ and every snapshot. Thus, $c(\tg)$ cops can play on $\tg'$ so that if at any time the robber moves to $P$ the cops act as if the robber moved to $u$. Since $u$ is a cutvertex between $P$ and $G$, eventually the robber is either captured on $u$, since $c(\tg) = c(\tg'[G])$, or it is somewhere on $P$ while a cop is on $u$. This latter cop eventually makes the catch. Therefore, $c(\tg')\leq c(\tg)$ and the equality holds. Similarly, $c(H) = c(G)$ and $c(H_{\max}) = c(\Sup)$. 
\qed \end{proof}

\subsection{Constructions based on the Petersen graph}

Some footprints are particularly useful and we derived a lot of constructions from the Petersen graph. We suspect that other graphs can help us fill the table. 

We employ the following graph operations multiple times in the next constructions. Given two graphs $H_1 = (V(H_1), E(H_1))$ and $H_2 = (V(H_2), E(H_2))$, the union of $H_1$ and $H_2$ is the graph $G = (V(H_1)\cup V(H_2), E(H_1)\cup E(H_2))$, written $G = H_1\cup H_2$. The \emph{join} of $H_1$ and $H_2$ is the graph $G = H_1\cup H_2 \cup \left(V(H_1)\cup V(H_2), \set{xy; x\in V(H_1), y\in V(H_2)}\right)$, written $G = H_1 + H_2$.

\begin{figure}
\centering
\includegraphics[scale=1]{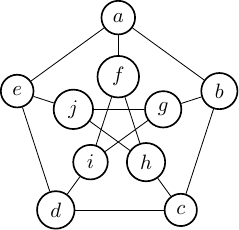}
\caption{A labelled Petersen graph}
\label{fig:pet_footprint}
\end{figure}

\begin{lemma}[(1,3,2)-copwin]\label{lem:132_copwin}
There exists a periodic graph $\tg$ with $c(G) < c(\tg) < c(\Sup)$.
\end{lemma}
\begin{proof}
Let $\pet$ be the Petersen graph, labelled as in \cref{fig:pet_footprint} and $x$ be a node that is not in $V(\pet)$. Let $O$ be the cycle $(a,b,c,d,e)$ and $I$ the cycle $(f,h,j,g,i)$. Let $G = \pet + (\set{x}, \emptyset)$. 
Let $H\subset G$ be the subgraph with $E(H) = E(O)\cup\set{\uv{a}{f}, \uv{b}{g}, \uv{c}{h}, \uv{d}{i}, \uv{e}{j}, \uv{a}{x}}$. Let $\tg$ be the following periodic graph. For every time $t\equiv 0\pmod{5}$, choose an edge $\uv{x}{u}\in E(G)$ and let $G_t = \pet \cup (\set{x,u}, \set{\uv{x}{u}}$). For every $t\not\equiv 0\pmod{5}$, let $G_t = H$. Observe that every snapshot is connected.

Since $x$ is a universal vertex in $G$, $c(G) = 1$. Moreover, $c(\Sup) = c(G_0) = 3$. Indeed, $c(G_t)\leq 3$ for every time $t$ since the smallest $4$-copwin graph has $19$ nodes \cite{turcotte20214} while $G$ has $11$ nodes. Furthermore, $G_0 = \pet \cup (\set{x,u},\set{xu})$ for some edge $\uv{x}{u}\in E(G)$, so $\pet$ is a retract of $G_0$ and $3 = c(\pet)\leq c(G_0)$ by the result of \citet{berarducci1993cop}. Finally, a single cop gains nothing in starting on $x$. Therefore, the robber can escape from a single cop by forever moving on $O$ since $c(O) = 2$ and $O \subset G_t$ for every time $t$. But, $O$ appears long enough in $\tg$ for two cops to make the catch. Therefore, $c(\tg) = c(O) = 2$. 
\qed \end{proof}

\begin{proposition}[(2,3,1)-copwin]\label{prop:231_copwin}
There exists a periodic graph $\tg$ with $c(\tg) < c(G) < c(\Sup)$.
\end{proposition}
\begin{proof}
Let $\pet$ be the Petersen graph and $x,y$ be two nodes that are not in $V(\pet)$. Let $O$ be the cycle $(a,b,c,d,e)$ and $I$ the cycle $(f,h,j,g,i)$. Then, let $G = \pet \cup (O + (\set{x}, \emptyset))\cup (I + (\set{y}, \emptyset))$. Observe that $\domnum(G) = 2$, with $\set{x,y}$ as the minimal dominating set. Therefore, $c(G) = 2$.

Let $\tg$ be the following periodic graph. For every time $t \equiv 0\pmod{11}$, let $e_t$ be an edge from $x$ to $O$, $f_t$ be an edge from $y$ to $I$ and $G_t = \pet\cup (V(G), \set{e_t,f_t})$. Thus, every such $G_t$ is connected. There exists a spanning tree $T$ of $\pet$ such that $T$ contains a node $u$ at distance at most $3$ from every other node in $T$. Let us add two edges connecting $x$ and $y$ to $T$ so that $u$ is at distance at most $4$ from either node. Let $G_t = T$ for every $t\not \equiv 0 \pmod{11}$. Then, $G_t$ is a spanning tree of $G$ that appears for at least $9>4$ consecutive snapshots, so a single cop can win on $\tg$. That is, $c(\tg) = 1$. Nevertheless, for every $t \equiv 0 \pmod{11}$, $\pet$ is a retract of $G_t$ given by mapping $x$ to its unique neighbour and similarly for $y$. Then, by Berarducci and Intrigila's classic result that $c(H)\leq c(H')$ whenever $H$ is a retract of $H'$ (\cite{berarducci1993cop}), $3 = c(\pet) \leq c(G_t)$. Finally, the smallest $4$-copwin graph has $19$ nodes (\cite{turcotte20214}), so $c(G_t)\leq 3$ and $c(G_t) = c(\Sup) = 3$. Therefore, $\tg$ is $(2,3,1)$-copwin.
\qed \end{proof}

\begin{figure}
\centering
\begin{subfigure}[ht]{0.3\linewidth}
\centering
\includegraphics[scale=1]{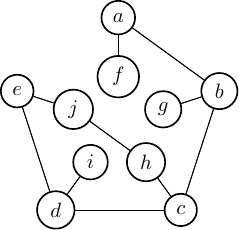}
\caption{$H_0$}
\end{subfigure}
\hfill
\begin{subfigure}[ht]{0.3\linewidth}
\centering
\includegraphics[scale=1]{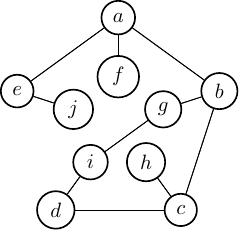}
\caption{$H_1$}
\end{subfigure}
\hfill
\begin{subfigure}[ht]{0.3\linewidth}
\centering
\includegraphics[scale=1]{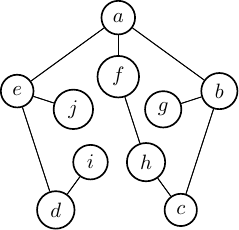}
\caption{$H_2$}
\end{subfigure}
\hfill
\begin{subfigure}[ht]{0.3\linewidth}
\centering
\includegraphics[scale=1]{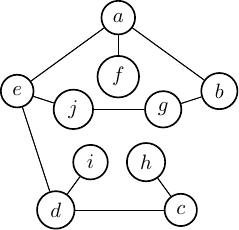}
\caption{$H_3$}
\end{subfigure}
\hfill
\begin{subfigure}[ht]{0.3\linewidth}
\centering
\includegraphics[scale=1]{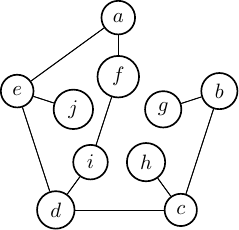}
\caption{$H_4$}
\end{subfigure}
\caption{Sequence of graphs used in \cref{lem:321_copwin}}
\label{fig:graphs_prop_321}
\end{figure}
\begin{lemma}[(3,2,1)-copwin]\label{lem:321_copwin}
There exists a periodic graph $\tg$ with $c(\tg) < c(\Sup) < c(G)$.
\end{lemma}
\begin{proof}
Let $G$ be the Petersen graph labelled as in \cref{fig:pet_footprint} and $H_0, \dots, H_4$ be the graphs shown in \cref{fig:graphs_prop_321}. Then, let $\tg = \peri{G_0,\dots,G_{p-1}}$ be the periodic graph where
\begin{align*}
G_t
&=
\begin{cases}
  H_0, &\text{ if } 0\leq t\leq 3\\
  H_1, &\text{ if } 4\leq t\leq 7\\
  H_2, &\text{ if } 8\leq t\leq 11\\
  H_3, &\text{ if } 12\leq t\leq 15\\
  H_4, &\text{ if } 16\leq t\leq 19.
\end{cases}
\end{align*}

The footprint is the Petersen graph, so $c(G) = 3$. Every snapshot is always connected and contains a cycle of length five, so $c(\Sup) \geq 2$. Moreover, two cops can capture the robber on every snapshot, so $c(\Sup)\leq 2$ and $c(\Sup) = 2$. 

Finally, consider the first $8$ snapshots $G_0,\dots,G_7$. Let a single cop start on $c$. The robber cannot start in $N_{H_0}[c] = \set{b,c,h,d}$. Neither can it start in $\set{a,f,g}$, since it would be stuck on a tree rooted at $c$. Finally, if it starts on $i$, it gets stuck on $i$ when the cop moves to $d$. Thus, the robber must start on $e$ or $j$. Let the cop wait until $G_4$, when $H_1$ appears. The robber cannot move out of the set of nodes $\set{e,j,h,c,d,i}$. In order for the robber to move to $i$ before $G_4$, it must move to $d$, in which case the cop makes the catch on the next snapshot. Similarly, $h\in \neig{4}{c}{\tg}$. Thus, the robber cannot safely end its turn on $\set{d,i,h}$ in $G_3$. Therefore, the robber must end its turn in $G_3$ somewhere in $\set{e,j}$. In this case, let the cop move to $b$ in $G_4$. The robber is on a tree rooted at $b$ until snapshot $G_7$ and the cop makes the catch.

This shows that $c(\tg) = 1$.
\qed \end{proof}

\subsection{Arguments from the non-existence of temporal corners}
\label{subsec:non_exist_temp_corners}

\begin{figure}
\centering
\begin{subfigure}[ht]{0.4\linewidth}
\includegraphics[scale=1]{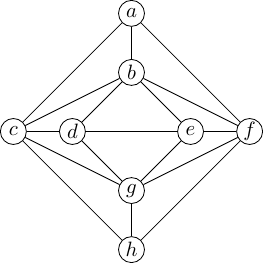}
\caption{Footprint $G$}
\end{subfigure}
\vfill
\begin{subfigure}[ht]{0.4\linewidth}
\includegraphics[scale=1]{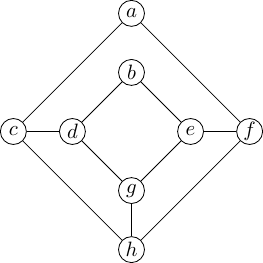}
\caption{Snapshot $G_0$}
\end{subfigure}
\hfill
\begin{subfigure}[ht]{0.4\linewidth}
\includegraphics[scale=1]{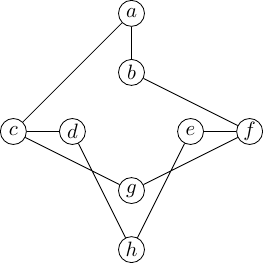}
\caption{Snapshot $G_1$}
\end{subfigure}
\vfill
\begin{subfigure}[ht]{0.4\linewidth}
\includegraphics[scale=1]{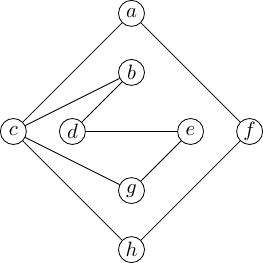}
\caption{Snapshot $G_2$}
\end{subfigure}
\caption{The footprint $G$ is copwin, the periodic graph $\tg = \peri{G_0,G_1,G_2}$ is not}
\label{fig:copwdyngraph}
\end{figure}

\begin{lemma}[(1,2,2)-copwin]\label{lem:copw_footprint_2sn_noncopwin}
There exists an always connected periodic graph with copwin footprint and $2$-copwin snapshots that is $2$-copwin.
\end{lemma}
\begin{proof}
The periodic graph $\tg$ in \cref{fig:copwdyngraph} is formed with $2$-copwin snapshots and copwin footprint. We can see every snapshot is $2$-copwin since each of them has girth $4$ and less than ten nodes. That is, none of them is copwin nor $3$-copwin, so all are $2$-copwin. 

At most two cops are necessary to capture the robber on $\tg$. Indeed, two cops can start on $d$ and $f$ in $G_0$, which is a dominating set of $G_0$. Then, $c(\tg)\leq \domnum(G_0) = 2$ because no matter where the robber starts, it gets captured in $G_0$. 

Inspection shows that $\tg$ has no temporal corner, so it cannot be copwin by \cref{obs:temp_corner}. Therefore, $2$ cops are necessary and sufficient to capture the robber.
\qed \end{proof}

\begin{figure}
\centering
\foreach \i in {0,...,4}{%
    \begin{subfigure}{0.3\textwidth}
    \centering
    \includegraphics[scale=0.7]{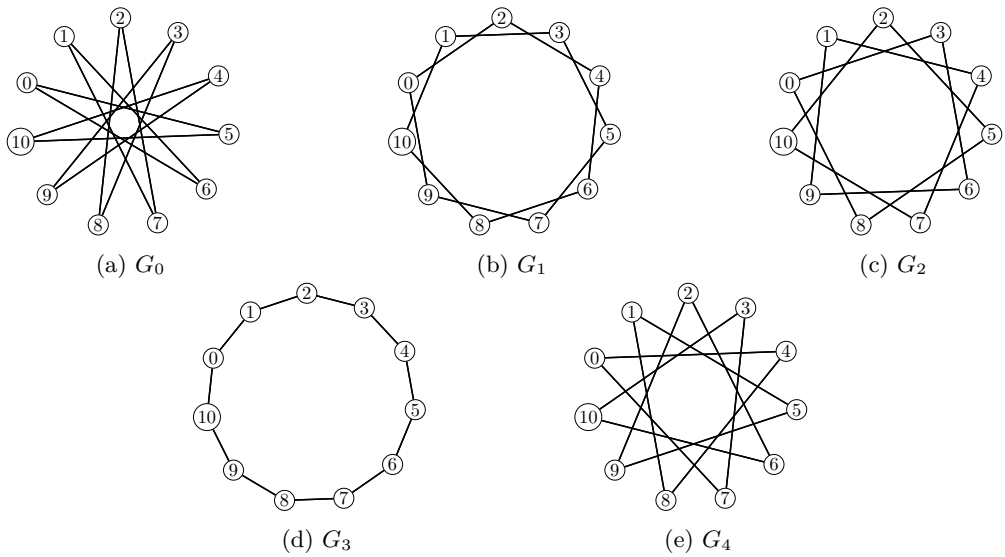}
    \hfill
    \protect\caption{$G_{\i}$}
    \end{subfigure}
}
\caption{The periodic graph used in \cref{prop:123-copwin_generic}. Each snapshot is a cycle of length $11$, the footprint is copwin and the periodic graph is $3$-copwin}
\label{fig:123_copwin}
\end{figure}

\begin{figure}
\centering
\includegraphics[scale=1]{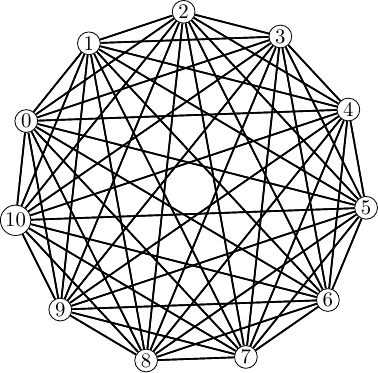}
\caption{The footprint $G$ of the periodic graph shown in \cref{fig:123_copwin} is a complete graph $K_{11}$}
\label{fig:123_copwin_footprint}
\end{figure}

\begin{proposition}[(1,2,3)-copwin]\label{prop:123-copwin_generic}
For any integers $n\geq 11$ and $p\geq 5$ odd, there exists an always connected periodic graph $\tg$ with period $p$ that is $(1,2,3)$-copwin.
\end{proposition}
\begin{proof}
Let $\tg$ be the periodic graph shown in \cref{fig:123_copwin}. One can verify that $\tg$ has no $2$-temporal corner since $\npg{t}{u}{\tg}$ is of the form $\set{u, u+i_t, u-i_t}$ such that $i_t \neq i_{t+1}$ for every time $t$. Thus, $c(\tg)\geq 3$ by \cref{prop:k_temp_corner}. Moreover, the footprint is shown in \cref{fig:123_copwin_footprint} and one can verify that every node has degree $10$ in the footprint. Thus, $G \equiv K_{11}$ is copwin. 

We can add more cycles $C_{11}$ to $\tg$ without creating $2$-temporal corners by adding pairs $(G_3,G_4)$ to the end of $\tg$. Note that this only generates odd periods. The \copnum{} of $G$ cannot go lower. Similarly, we can always apply \cref{prop:abccop_larger_abccop} to increase the number of nodes without changing the properties of $\tg$. 

Let us show that $c(\tg)\leq 3$. Let us place one cop on node $0$, one on node $3$ and another on $8$. Let the cops wait until $G_3$. In $G_3$, the robber must be on either $5$ or $6$ not to be captured in $G_3$. Then, $\neig{3}{5}{\tg}\cup \neig{3}{6}{\tg} = \set{4,5,6,7} \subseteq \neig{4}{0}{\tg} \cup \neig{4}{2}{\tg} \cup \neig{4}{9}{\tg} = \set{0,2,4,5,6,7,9}$, so when the cops in $G_3$ move to $0$, $2$ and $9$ the robber gets stuck on a $3$-temporal corner of the cops' positions. Then, the cops win in $G_4$.
\qed \end{proof}

The rest of the results are stated in the following theorem, whose proof can be found in \citet{simard2024temporal}.
\begin{theorem}\label{prop:leftover_results}
    There exist periodic graphs that are $(a,b,c)$-copwin with:
    \begin{align*}
    (a,b,c) &\in
    \set*{
        \begin{array}{cccc}
        (1,1,1) & (1,2,1) & (1,3,1) & (1,3,3)\\
        (2,1,1) & (2,1,2) & (2,2,2) & (2,2,3)\\
        (2,3,2) & (2,3,3) & (3,1,2) & (3,2,2)\\ 
        (3,3,1) & (3,2,3) & (3,3,2) & (3,3,3)
        \end{array}
    }
    \end{align*}
\end{theorem}

\section{Conclusion}

We showed that allowing a graph to change over time generates new challenges when studying the game of Cops and Robber. We presented many results to inform researchers of those challenges. We also exhibited results on periodic graphs that generalize their counterparts on static graphs, which shows what type of extra assumptions can be used to extend them.

This work can be seen as a first foray into the analysis of the game of Cops and Robber on Periodic Graphs. One major avenue of research this opens up is the determination of \copnum{s} of specific classes of \emph{periodic graphs}. This is common on undirected graphs. Recall for example that we know the \copnum{s} of planar graphs \cite{Aigner1984}, outerplanar graphs \cite{clarke2002constrained} and hypercubes \cite{maamoun1987game}. We gave partial answers to this question by focusing on $\ctgmax{G}$: then, for example, every periodic graph with footprint $G$ has \copnum{} at most $\tw(G) + 1$ (\cref{thm:treewidth_foot_pg}). In order to answer this question on classes of periodic graphs that are not purely defined by their footprints, one would have to come up with interesting classes of periodic graphs. This is a major hurdle we faced. The literature so far is sparse on classes of periodic graphs. Moreover, those classes that have been defined, such as those presented by Casteigts et al. \cite{Casteigts2018}, are often not well suited to the game of Cops and Robber so it is difficult to compute their \copnum{s}. 

The class of periodic graphs with footprint $G$, which we investigated when computing $\ctgmax{G}$, is nevertheless vast and interesting. One confounding aspect of this class of periodic graphs is that, although a priori it looks like computing $\ctgmax{G}$ is akin to playing the game with imperfect information, whenever an element of this family is chosen the cops will know of it. This is similar to computing the maximum \copnum{} of a family of random graphs: the structure is only known a posteriori. Furthermore, the value of $\ctgmax{G}$ says something about the nature of $G$. The fact that when $G$ is outerplanar we have $\ctgmax{G} - c(G) \leq 1$ means that outerplanar graphs have particularly strong separation properties. To emphasize: when $G$ is outerplanar, no matter what periodic sequence of subgraphs of $G$ one takes, $3$ cops can always capture the robber. This is in line with the type of strategies we used to prove that $\ctgmax{G}\leq 3$. One might think of those strategies as \emph{stubborn}: they tell the cops where to go on $G$ and the cops make their moves as their incident edges become available. The $3$-copwin strategy on planar graphs, recall \cite{Aigner1984}, are \emph{dynamic} because they involve moving the cops, possibly at every turn, to guard isometric paths. We suspect that $\ctgmax{G} - c(G)$ grows larger with the number of vertices of $G$ when $G$ is planar. This suggests the following line of reasoning. So far, when computing $\ctgmax{G}$ we have been looking at properties of $G$ to understand properties of the class of periodic graphs with footprint $G$. However, one could also seek to understand properties of $G$ from properties of the periodic graphs with footprint $G$. For example, if it turns out that $\ctgmax{G}$ is not bounded by a constant when $G$ is planar, then one might conclude that the \emph{only} winning strategies for the cops on $G$ are dynamic.

\bibliographystyle{abbrvnat}
\bibliography{pgcr_tidied}

\begin{thebibliography}{53}
\providecommand{\natexlab}[1]{#1}
\providecommand{\url}[1]{\texttt{#1}}
\expandafter\ifx\csname urlstyle\endcsname\relax
  \providecommand{\doi}[1]{doi: #1}\else
  \providecommand{\doi}{doi: \begingroup \urlstyle{rm}\Url}\fi

\bibitem[Agarwalla et~al.(2018)Agarwalla, Augustine, Moses~Jr, Madhav, and
  Sridhar]{AgMMS18}
A.~Agarwalla, J.~Augustine, W.~K. Moses~Jr, S.~K. Madhav, and A.~K. Sridhar.
\newblock Deterministic dispersion of mobile robots in dynamic rings.
\newblock In \emph{19th International Conference on Distributed Computing and
  Networking}, pages 1--4, 2018.

\bibitem[Aigner and Fromme(1984)]{Aigner1984}
M.~Aigner and M.~Fromme.
\newblock {A game of cops and robbers}.
\newblock \emph{Discrete Applied Mathematics}, 8\penalty0 (1):\penalty0 1--12,
  1984.

\bibitem[Baird et~al.(2014)Baird, Beveridge, Bonato, Codenotti, Maurer,
  McCauley, and Valeva]{baird2013minimum}
W.~Baird, A.~Beveridge, A.~Bonato, P.~Codenotti, A.~Maurer, J.~McCauley, and
  S.~Valeva.
\newblock {On the minimum order of k-cop win graphs}.
\newblock \emph{Contributions to Discrete Mathematics}, 9\penalty0 (1), 2014.

\bibitem[Balev et~al.(2020)Balev, Laredo, Lamprou, Pigné, and
  Sanlaville]{Balev2020}
S.~Balev, J.~L.~J. Laredo, I.~Lamprou, Y.~Pigné, and E.~Sanlaville.
\newblock Cops and robbers on dynamic graphs: Offline and online case.
\newblock In \emph{27th International Colloquium on Structural Information and
  Communication Complexity (SIROCCO)}, pages 203–--219, 2020.

\bibitem[Berarducci and Intrigila(1993)]{berarducci1993cop}
A.~Berarducci and B.~Intrigila.
\newblock {On the cop number of a graph}.
\newblock \emph{Advances in Applied Mathematics}, 14\penalty0 (4):\penalty0
  389--403, 1993.

\bibitem[Berwanger(2013)]{Berwanger2013GraphGW}
D.~Berwanger.
\newblock {\em {Graph games with perfect information}}.
\newblock Lecture notes, ENS Paris-Saclay, 2013.

\bibitem[Bodlaender(1996)]{bodlaender1996linear}
H.~L. Bodlaender.
\newblock {A linear-time algorithm for finding tree-decompositions of small
  treewidth}.
\newblock \emph{Journal on Computing}, 25\penalty0 (6):\penalty0 1305--1317,
  1996.

\bibitem[Bonato and Mac{G}illivray(2017)]{BoM17}
A.~Bonato and G.~Mac{G}illivray.
\newblock {Characterizations and algorithms for generalized Cops and Robbers
  games}.
\newblock \emph{Contributions to Discrete Mathematics}, 12\penalty0 (1), 2017.

\bibitem[Bonato and Nowakowski(2011)]{Bonato2011f}
A.~Bonato and R.~Nowakowski.
\newblock \emph{{The Game of Cops and Robbers on Graphs}}.
\newblock American Mathematical Society, 2011.

\bibitem[Bournat et~al.(2017)Bournat, Dubois, and Petit]{BouDP17}
M.~Bournat, S.~Dubois, and F.~Petit.
\newblock Computability of perpetual exploration in highly dynamic rings.
\newblock In \emph{IEEE 37th International Conference on Distributed Computing
  Systems (ICDCS)}, pages 794--804, 2017.

\bibitem[Casteigts(2018)]{Casteigts2018}
A.~Casteigts.
\newblock {\em {A Journey through Dynamic Networks (with Excursions)}}.
\newblock Habilitation à diriger des recherches, Université de Bordeaux,
  2018.

\bibitem[Casteigts et~al.(2012)Casteigts, Flocchini, Quattrociocchi, and
  Santoro]{Casteigts2013}
A.~Casteigts, P.~Flocchini, W.~Quattrociocchi, and N.~Santoro.
\newblock {Time-varying graphs and dynamic networks}.
\newblock \emph{International Journal of Parallel, Emergent and Distributed
  Systems}, 27\penalty0 (5):\penalty0 387--408, 2012.

\bibitem[Casteigts et~al.(2014)Casteigts, Flocchini, Mans, and
  Santoro]{CaFMS14}
A.~Casteigts, P.~Flocchini, B.~Mans, and N.~Santoro.
\newblock Measuring temporal lags in delay-tolerant networks.
\newblock \emph{IEEE Transactions on Computers}, 63\penalty0 (2):\penalty0
  397--410, 2014.

\bibitem[Clarke(2002)]{clarke2002constrained}
N.~E. Clarke.
\newblock \emph{{Constrained cops and robber}}.
\newblock PhD thesis, Dalhousie University, 2002.

\bibitem[Clarke and Nowakowski(2001)]{Clarke2001}
N.~E. Clarke and R.~J. Nowakowski.
\newblock {Cops, robber and traps}.
\newblock \emph{Utilitas Mathematica}, 60:\penalty0 91--98, 2001.

\bibitem[De~Carufel et~al.(2023)De~Carufel, Flocchini, Santoro, and
  Simard]{CaFlSaSi23}
J.-L. De~Carufel, P.~Flocchini, N.~Santoro, and F.~Simard.
\newblock Cops \& robber on periodic temporal graphs: Characterization and
  improved bounds.
\newblock In \emph{30th International Conference on Structural Information and
  Communication Complexity (SIROCCO)}, pages 29--41, 2023.

\bibitem[Di~Luna(2019)]{DiL19}
G.~A. Di~Luna.
\newblock Mobile agents on dynamic graphs.
\newblock In P.~Flocchini, G.~Prencipe, and N.~Santoro, editors,
  \emph{Distributed Computing by Mobile Entities}, pages 549--584. Springer
  International Publishing, 2019.

\bibitem[Di~Luna et~al.(2017)Di~Luna, Flocchini, Pagli, Prencipe, Santoro, and
  Viglietta]{DiLFPPSV17}
G.~A. Di~Luna, P.~Flocchini, L.~Pagli, G.~Prencipe, N.~Santoro, and
  G.~Viglietta.
\newblock Gathering in dynamic rings.
\newblock In \emph{24th International Conference on Structural Information and
  Communication Complexity (SIROCCO)}, pages 339--355, 2017.

\bibitem[Di~Luna et~al.(2020)Di~Luna, Dobrev, Flocchini, and Santoro]{DiLDFS20}
G.~A. Di~Luna, S.~Dobrev, P.~Flocchini, and N.~Santoro.
\newblock {Distributed exploration of dynamic rings}.
\newblock \emph{Distributed Computing}, 33:\penalty0 41--67, 2020.

\bibitem[Diestel(2000)]{Diestel2000}
R.~Diestel.
\newblock \emph{{Graph Theory}}.
\newblock Springer Berlin, Heidelberg, 2000.

\bibitem[Erlebach and Spooner(2018)]{ErS18}
T.~Erlebach and J.~T. Spooner.
\newblock Faster exploration of degree-bounded temporal graphs.
\newblock In \emph{{43rd International Symposium on Mathematical Foundations of
  Computer Science (MFCS)}}, pages 1--13, 2018.

\bibitem[Erlebach and Spooner(2020)]{ErS20}
T.~Erlebach and J.~T. Spooner.
\newblock A game of cops and robbers on graphs with periodic edge-connectivity.
\newblock In \emph{46th International Conference on Current Trends in Theory
  and Practice of Informatics (SOFSEM)}, pages 64--75, 2020.

\bibitem[Erlebach et~al.(2021)Erlebach, Hoffmann, and
  Kammer]{erlebach2021temporal}
T.~Erlebach, M.~Hoffmann, and F.~Kammer.
\newblock {On temporal graph exploration}.
\newblock \emph{Journal of Computer and System Sciences}, 119:\penalty0 1--18,
  2021.

\bibitem[Erlebach et~al.(2024)Erlebach, Morawietz, Spooner, and
  Wolf]{erlebach2024}
T.~Erlebach, N.~Morawietz, J.~T. Spooner, and P.~Wolf.
\newblock A cop and robber game on edge-periodic temporal graphs.
\newblock \emph{Journal of Computer and System Sciences}, 144:\penalty0 103534,
  2024.

\bibitem[Ferreira(2004)]{Fe04}
A.~Ferreira.
\newblock {Building a reference combinatorial model for MANETs}.
\newblock \emph{{IEEE Network}}, 18\penalty0 (5):\penalty0 24--29, 2004.

\bibitem[Flocchini et~al.(2013)Flocchini, Mans, and
  Santoro]{flocchini2013exploration}
P.~Flocchini, B.~Mans, and N.~Santoro.
\newblock {On the exploration of time-varying networks}.
\newblock \emph{Theoretical Computer Science}, 469:\penalty0 53--68, 2013.

\bibitem[Fluschnik et~al.(2020)Fluschnik, Molter, Niedermeier, Renken, and
  Zschoche]{fluschnik2020time}
T.~Fluschnik, H.~Molter, R.~Niedermeier, M.~Renken, and P.~Zschoche.
\newblock As time goes by: {R}eflections on treewidth for temporal graphs.
\newblock In F.~V. "Fomin, S.~Kratsch, and E.~J. van Leeuwen, editors,
  \emph{Treewidth, Kernels, and Algorithms}, pages 49--77. Springer
  International Publishing, 2020.

\bibitem[Fomin et~al.(2010)Fomin, Golovach, Kratochvil, Nisse, and
  Suchan]{fomin2010pursuing}
F.~V. Fomin, P.~A. Golovach, J.~Kratochvil, N.~Nisse, and K.~Suchan.
\newblock {Pursuing a fast robber on a graph}.
\newblock \emph{Theoretical Computer Science}, 411\penalty0 (7-9):\penalty0
  1167--1181, 2010.

\bibitem[Gotoh et~al.(2018)Gotoh, Sudo, Ooshita, Kakugawa, and
  Masuzawa]{GoSOKM18}
T.~Gotoh, Y.~Sudo, F.~Ooshita, H.~Kakugawa, and T.~Masuzawa.
\newblock Group exploration of dynamic tori.
\newblock In \emph{IEEE 38th International Conference on Distributed Computing
  Systems (ICDCS)}, pages 775--785, 2018.

\bibitem[Gotoh et~al.(2021)Gotoh, Flocchini, Masuzawa, and Santoro]{GoSOKM21}
T.~Gotoh, P.~Flocchini, T.~Masuzawa, and N.~Santoro.
\newblock Exploration of dynamic networks: tight bounds on the number of
  agents.
\newblock \emph{Journal of Computer and System Sciences}, 122:\penalty0 1--18,
  2021.

\bibitem[Harary and Gupta(1997)]{HarG97}
F.~Harary and G.~Gupta.
\newblock Dynamic graph models.
\newblock \emph{Mathematical and Computer Modelling}, 25\penalty0 (7):\penalty0
  79--88, 1997.

\bibitem[Harvey and Wood(2017)]{harvey2017parameters}
D.~J. Harvey and D.~R. Wood.
\newblock {Parameters tied to treewidth}.
\newblock \emph{Journal of Graph Theory}, 84\penalty0 (4):\penalty0 364--385,
  2017.

\bibitem[Hill(2008)]{hill2008cops}
A.~Hill.
\newblock \emph{{Cops and robbers: theme and variations}}.
\newblock PhD thesis, Dalhousie University, 2008.

\bibitem[Holme and Saramäki(2012)]{Holme2012}
P.~Holme and J.~Saramäki.
\newblock {Temporal networks}.
\newblock \emph{Physics Reports}, 519\penalty0 (3):\penalty0 97--125, 2012.

\bibitem[Ilcinkas and Wade(2011)]{IlW11}
D.~Ilcinkas and A.~M. Wade.
\newblock On the power of waiting when exploring public transportation systems.
\newblock In \emph{15th International Conference on Principles of Distributed
  Systems (OPODIS)}, pages 451--464, 2011.

\bibitem[Ilcinkas et~al.(2014)Ilcinkas, Klasing, and Wade]{IlKW14}
D.~Ilcinkas, R.~Klasing, and A.~M. Wade.
\newblock Exploration of constantly connected dynamic graphs based on cactuses.
\newblock In \emph{21st International Colloquium on Structural Information and
  Communication Complexity (SIROCCO)}, pages 250--262, 2014.

\bibitem[Jathar et~al.(2014)Jathar, Yadav, and Gupta]{JatYG14}
R.~Jathar, V.~Yadav, and A.~Gupta.
\newblock {Using periodic contacts for efficient routing in delay tolerant
  networks}.
\newblock \emph{Ad Hoc \& Sensor Wireless Networks}, 22\penalty0
  (1,2):\penalty0 283--308, 2014.

\bibitem[Joret et~al.(2010)Joret, Kamiński, and Theis]{joret2008cops}
G.~Joret, M.~Kamiński, and D.~O. Theis.
\newblock {The cops and robber game on graphs with forbidden (induced)
  subgraphs}.
\newblock \emph{Contributions to Discrete Mathematics}, 5\penalty0 (2), 2010.

\bibitem[Khatri et~al.(2018)Khatri, Komarov, Krim-Yee, Kumar, Seamone, Virgile,
  and Xu]{Khatri2018}
D.~Khatri, N.~Komarov, A.~Krim-Yee, N.~Kumar, B.~Seamone, V.~Virgile, and
  A.~Xu.
\newblock {A study of cops and robbers in oriented graphs}.
\newblock \emph{arXiv preprint, arXiv:1811.06155}, pages 1--23, 2018.

\bibitem[Kinnersley(2015-03)]{Kinnersley2015}
W.~B. Kinnersley.
\newblock Cops and robbers is exptime-complete.
\newblock \emph{Journal of Combinatorial Theory, Series B}, 111:\penalty0
  201--220, 2015-03.

\bibitem[Kuhn et~al.(2010)Kuhn, Lynch, and Oshman]{KuhLyO10}
F.~Kuhn, N.~Lynch, and R.~Oshman.
\newblock {Distributed computation in dynamic networks}.
\newblock In \emph{42nd ACM Symposium on Theory of Computing}, pages 513--522,
  2010.

\bibitem[Loh and Oh(2017)]{loh2017cops}
P.-S. Loh and S.~Oh.
\newblock {Cops and Robbers on Planar-Directed Graphs}.
\newblock \emph{Journal of Graph Theory}, 86\penalty0 (3):\penalty0 329--340,
  2017.

\bibitem[Maamoun and Meyniel(1987)]{maamoun1987game}
M.~Maamoun and H.~Meyniel.
\newblock On a game of policemen and robber.
\newblock \emph{Discrete Applied Mathematics}, 17\penalty0 (3):\penalty0
  307--309, 1987.

\bibitem[Mehrabian(2011)]{mehrabian2011cops}
A.~Mehrabian.
\newblock {\em {Cops and robber game with a fast robber}}.
\newblock Master's thesis, University of Waterloo, 2011.

\bibitem[Morawietz and Wolf(2021)]{morawietz2021timecop}
N.~Morawietz and P.~Wolf.
\newblock A timecop's chase around the table.
\newblock \emph{arXiv preprint, arXiv:2104.08616}, 2021.

\bibitem[Morawietz et~al.(2020)Morawietz, Rehs, and Weller]{morawietz2020}
N.~Morawietz, C.~Rehs, and M.~Weller.
\newblock {A timecop's work is harder than you think}.
\newblock In \emph{45th International Symposium on Mathematical Foundations of
  Computer Science (MFCS)}, pages 71.1--71.14, 2020.

\bibitem[Nowakowski and Winkler(1983)]{Nowakowski1983}
R.~J. Nowakowski and P.~Winkler.
\newblock {Vertex-to-vertex pursuit in a graph}.
\newblock \emph{Discrete Mathematics}, 43\penalty0 (2-3):\penalty0 235--239,
  1983.

\bibitem[O'Dell and Wattenhofer(2005)]{OdW05}
R.~O'Dell and R.~Wattenhofer.
\newblock {Information dissemination in highly dynamic graphs}.
\newblock In \emph{Joint Workshop on Foundations of Mobile Computing}, pages
  104--110, 2005.

\bibitem[Quilliot(1978)]{Qu78}
A.~Quilliot.
\newblock \emph{{Jeux et points fixes sur les graphes}}.
\newblock Thèse de 3ème cycle, Université de Paris VI, 1978.

\bibitem[Simard(2024)]{simard2024temporal}
F.~Simard.
\newblock \emph{Temporal distances and cops and robber games on temporal
  networks}.
\newblock PhD thesis, Université d'Ottawa $\vert$ University of Ottawa, 2024.

\bibitem[Simard et~al.(2021)Simard, Desharnais, and
  Laviolette]{simard2021general}
F.~Simard, J.~Desharnais, and F.~Laviolette.
\newblock General cops and robbers games with randomness.
\newblock \emph{Theoretical Computer Science}, 887:\penalty0 30--50, 2021.

\bibitem[Turcotte and Yvon(2021)]{turcotte20214}
J.~Turcotte and S.~Yvon.
\newblock {4-cop-win graphs have at least 19 vertices}.
\newblock \emph{Discrete Applied Mathematics}, 301:\penalty0 74--98, 2021.

\bibitem[Wehmuth et~al.(2015)Wehmuth, Ziviani, and Fleury]{WeZF15}
K.~Wehmuth, A.~Ziviani, and E.~Fleury.
\newblock {A unifying model for representing time-varying graphs}.
\newblock In \emph{2nd IEEE International Conference on Data Science and
  Advanced Analytics, (DSAA)}, pages 1--10, 2015.

\end{thebibliography}



\end{document}